\documentclass[11pt, reqno]{amsart}
\pdfoutput=1
\usepackage{amssymb,amsmath,amsthm,amsfonts, upgreek,latexsym}
\usepackage{epsfig,esint}
\usepackage{mathrsfs}
\usepackage{mathtools,textcomp}
\usepackage{graphicx}
\usepackage{epic,eepic,wrapfig,color, ifthen} 
\usepackage{url,cite}
\usepackage{systeme}
 
\DeclareFontFamily{U}{mathx}{}
\DeclareFontShape{U}{mathx}{m}{n}{<-> mathx10}{}
\DeclareSymbolFont{mathx}{U}{mathx}{m}{n}
\DeclareMathAccent{\wc}{0}{mathx}{"71}

\usepackage[colorlinks=true, citecolor=blue, filecolor=cyan, linkcolor=red, urlcolor=magenta]{hyperref}
\usepackage{pmboxdraw}
\usepackage{verbatim}
\usepackage[normalem]{ulem}
\usepackage{enumerate}
\mathtoolsset{showonlyrefs}
\usepackage[bottom]{footmisc}
\newtheorem{thm}{Theorem}[section]

\newtheorem{lemma}[thm]{Lemma}

\newtheorem{prop}[thm]{Proposition}
\newtheorem{corollary}[thm]{Corollary}
\newtheorem{defn}[thm]{Definition}
\newtheorem{example}[thm]{Example}

\def\diam{{\text{\rm diam}}}

\numberwithin{equation}{section}
\allowdisplaybreaks

\newcommand{\bk}{\color{black}}

\def\bn{{\bf n}}

\def\1{{\bf 1}}
\def\ind{{\bf 1}}
\def\nn{\nonumber}
\def\bee{\begin{equation}}
\def\eee{\end{equation}}
\def\sA {{\mathcal A}} \def\sB {{\mathcal B}} 
 \def\sE {{\mathcal E}} \def\sF {{\mathcal F}}

 \def\sN {{\mathcal N}}

 \def\bH {{\mathbb H}}

\def\R {{\mathbb R}}  \def\bH {{\mathbb H}}
 \def\Z {{\mathbb Z}}

\newcommand{\overbar}[1]{\mkern 0.5mu\overline{\mkern-0.5mu#1\mkern-0.5mu}\mkern 0.5mu}

\def\eps{\varepsilon}
\def\e{\mathbf e}

\def\FF{{\mathcal F}}
\def\EE{{\mathcal E}}

\def\R{{\mathbb R}}

\def\E{{\mathbb E}}
\def\F{{\bf F}}
\def\P{{\mathbb P}}
\def\N{{\mathbb N}}
\def\eps{\varepsilon}
\def\wh{\widehat}
\def\wt{\widetilde}

\def\ub{{\overbar \beta}}
\begin{document}
	\title[Heat kernel estimates for degenerate Markov processes]
{ \bf   
Dirichlet heat kernel estimates for Markov processes with jump kernels blowing up at the boundary}

\author{Soobin Cho, \quad Panki Kim, \quad Renming Song \quad and \quad Zoran Vondra\v{c}ek}

\address[Cho]{School of Mathematics and Statistics, Beijing Institute of Technology, Beijing 100081, China}
\curraddr{}
\email{soobincho.math@gmail.com}

\address[Kim]{Department of Mathematical Sciences and Research Institute of Mathematics,
	Seoul National University,	Seoul 08826, Republic of Korea}
\curraddr{}
\thanks{Panki Kim is supported by the National Research Foundation of Korea(NRF) grant funded by the Korea government(MSIP) (No. RS-2023-00270314).}
\email{pkim@snu.ac.kr}

\address[Song]{
	Department of Mathematics, University of Illinois Urbana-Champaign, Urbana, IL 61801,
	USA}
\curraddr{}
\thanks{Research of Renming Song is supported in part by a grant from
	the Simons Foundation \#960480.}
\email{rsong@illinois.edu}

\address[Vondra\v{c}ek]
{
	Dr.~Franjo Tu{\dj}man Defense and Security 
	University
	 and Department of Mathematics, Faculty of Science, University of Zagreb, Zagreb, Croatia
}
\curraddr{}
\thanks{Zoran Vondra\v{c}ek is supported in part by the Croatian Science Foundation under the project IP-2025-02-8793 and by the European Union – NextGenerationEU through the National Recovery and Resilience Plan 2021–2026 Institutional grants of University of Zagreb Faculty of Science (PMF-PRESTIGE)} 
\email{vondra@math.hr}

\date{}

\begin{abstract}
Let $D$ be a $C^{1, \eta}$ open set in $\mathbb{R}^d$ and $\mathcal{B}: D\times D\to (0, \infty)$ a symmetric Borel  function which blows up at the boundary of $D$. Let $\overline{X}$ be a conservative purely discontinuous Markov process on $\overline{D}$ with jump kernel $\mathcal{B}(x,y)|x-y|^{-d-\alpha}$, $0<\alpha<2$. In a recent paper, we obtained two-sided estimates on the transition density of $\overline{X}$. In this paper, we consider the process $X$ obtained by killing  $\overline{X}$ upon 
approaching the boundary of $D$ or by a critical potential. We obtain sharp two-sided estimates on the transition density of $X$. The estimates capture the interplay between the blow-up rate of $\mathcal{B}$, the geometry of $D$, and killing.
\end{abstract}

\bigskip

\maketitle



\bigskip


\section{Introduction and main result}\label{s:intro}

The study of heat kernels for discontinuous Markov processes and the corresponding nonlocal operators has been a major area of research in probability theory and partial differential equations over the last few decades. Beginning with the pioneering work \cite{CK03} on stable-like processes, sharp heat kernel estimates have been established for a wide range of jump processes, in settings ranging from metric measure spaces to subsets of $\mathbb{R}^d$ with smooth or rough boundaries; see, for example, 
\cite{BGR10, BGR14, CKS10, CKS10a, CKW20, CKW21, GHL17, GHH18} 
and the references therein. A common assumption underlying much of this theory is a uniform bound on the tails of the associated jump measures.

Important examples of jump processes for which the uniform tail bound fails arise from the trace of a stable process on a smooth open set and from the process introduced in \cite{DRV17, Von21} in connection with the nonlocal Neumann problem. 
In both cases, the corresponding jump kernels blow up near the boundary, leading to unbounded tails of the associated jump measures.
The nonlocal Neumann problem has also been investigated from an analytic perspective, with various extensions of the nonlocal Neumann problem developed in \cite{Ab20, AFR20, MPL, DRVW26}.

These examples have prompted the development of a potential theory for purely discontinuous Markov processes in open sets $D\subset\mathbb{R}^d$ whose jump kernels have the form
$$
J(x,y)=\sB(x,y)|x-y|^{-d-\alpha},
$$
where $\sB(x,y)$ may blow up near the boundary of $D$. The first results in this direction were obtained in \cite{KSV23}, where potential-theoretic properties away from the boundary were investigated. The paper \cite{KSV25} studied the case when $D$ is the half-space and established the boundary Harnack principle and sharp Green function estimates. A substantially broader framework was subsequently developed in \cite{CK25+}. In that work, the stability of Hölder regularity for both caloric and harmonic functions was established in a fairly general setting. The authors also introduced a new formulation of the tail condition for jump measures and established corresponding weighted functional inequalities.

Building on these developments, in \cite{CKSV25} we established sharp two-sided heat kernel estimates for a broad class of purely discontinuous, conservative symmetric Markov 
processes $\overline X$ 
on closed subsets $F$ of $\mathbb{R}^d$. In this setting, the jump kernel may become unbounded along a Borel subset $\Sigma\subset F$ whose Assouad dimension is strictly less than $d$, while the complement $F\setminus\Sigma$ is assumed to be dense in $F$ and $\kappa$-fat. The analysis is complicated by the absence of a uniform bound on the tails of the jump measures. In this setting, the parabolic Harnack principle fails, and standard heat kernel methods do not provide the a priori off-diagonal estimates needed for the analysis.

In the present paper, we study heat kernel estimates for two subprocesses of the conservative process $\overline X$ in a $C^{1,\eta}$ open set $D\subset\mathbb{R}^d$: the process obtained by killing $\overline X$ upon approaching the boundary $\partial D$, and the process obtained by killing $\overline X$ with a critical killing potential. The corresponding elliptic theory was developed in the companion paper \cite{CKSV26a}, under certain additional assumptions that will be specified below.

The boundary blow-up of the jump intensity considered above is one way in which uniform ellipticity can fail. Another possibility is that the factor $\sB(x,y)$ tends to zero near the boundary. Sharp two-sided heat kernel estimates for purely discontinuous Markov processes with this type of boundary degeneracy in a $C^{1,1}$ open set $D\subset\mathbb{R}^d$ were recently established in \cite{CKSV25b}.

\medskip

We now describe the precise assumptions on the jump kernel, beginning with the blow-up mechanism that is central to the present paper. 
We say that $\Phi:[0,\infty) \to [1,\infty)$ is a \emph{blow-up weight function} if it is continuous, increasing, 
satisfies  $\Phi(0)=1$ and has the following upper scaling property: there exist constants $\beta\ge 0$ and $C>1$ such that
\begin{align}\label{e:upper-scaling}
	\Phi(r)/\Phi(s) \le C(r/s)^{ \beta} \quad \text{for all $0<s\le r$}.
\end{align}
The upper  Matuszewska index $\overline \beta$ of $\Phi$ is defined as  the infimum of all $\beta$ 
for which \eqref{e:upper-scaling} holds.

Let ${\rm dim}_{\rm A}(\partial D)$ be the Assouad dimension of $\partial D$. For the definition of the Assouad dimension, see \cite[Definition 2.1]{CKSV25}.
Suppose that $D\subset \R^d$ is a Lipschitz open set with localization radius $R_0$, and that  ${\rm dim}_{\rm A}(\partial D)$ is $d-1$.

Let
$\sB:D \times D \to (0,\infty)$ be a Borel function satisfying the following assumption:

\medskip

\noindent  \hypertarget{A}{{\bf (A)}} 
 There exist a blow-up weight function $\Phi$  with upper Matuszewska index $\overline \beta$ satisfying
\begin{align*}
\overline \beta<1 \wedge \alpha
\end{align*} 
and constants $A_0 \in (0,\infty]$ and $C_1>1$ such that for all $x,y \in D$,
\begin{align}\label{e:J-Phi}
	C_1^{-1}\Phi\left( \frac{(|x-y|\wedge A_0)^2}{(\delta_D(x)\wedge A_0) (\delta_D(y) \wedge A_0)}\right)\le 	\sB(x,y) 
	\le C_1\Phi\left( \frac{(|x-y|\wedge A_0)^2}{(\delta_D(x)\wedge A_0) (\delta_D(y) \wedge A_0)}\right),
\end{align}
where $\delta_D(x)$ stands for the distance between $x$ and $D^c$, and 
$a\wedge b:=\min \{a, b\}$ and  $a\vee b:=\max\{a, b\}$.
\bigskip

It follows that $\sB(x,x)\asymp 1$ and that $\sB(x,y)\ge C_1^{-1}$. 
Here, and throughout this paper, $A\asymp B$ means that 
 the quotient of $A$ and $B$
 is bounded between two positive constants.
We also remark that, in many examples, one can take the constant $A_0$ to be $\mathrm{diam}(D^c)$.

Consider the bilinear form
\begin{align*}
	\sE^0(u,v):= \frac12 \int_{D\times D} (u(x)-u(y))(v(x)-v(y)) \frac{\sB(x,y)}{|x-y|^{d+\alpha}} dxdy.
\end{align*}
This form was studied in \cite{CKSV25}, under weaker assumptions on the set $D$.

Let ${\rm Lip}_c(\overline D)$ be the family of all Lipschitz functions on $\overline D$ with compact supports. Since the $d$-dimensional Lebesgue measure of $\partial D$ is zero, we can identify $L^2(\overline D)$
with $L^2( D)$. Under \hyperlink{A}{{\bf (A)}}, $\sE^0(u,u)<\infty$ for all $u\in {\rm Lip}_c(\overline D)$; see \cite[Proposition 3.1]{CKSV25}.

Let  $\overline \sF$ be the closure of $\mathrm{Lip}_c(\overline D)$ in $L^2(D)$ under  the norm $(\sE^0_1)^{1/2}$, 
where  $\sE^0_1:=\sE^0 + \lVert \cdot \rVert_{L^2(D)}^2$.  Then  $(\sE^0, \overline \sF)$ is a regular Dirichlet 
form on $L^2(\overline D)$.
Let $\overline{X}=(\overline{X}_t)_{t\ge 0}$ be the conservative Hunt process associated with the form $(\sE^0, \overline \sF)$.
Joint H\"older continuity of bounded caloric functions (parabolic H\"older regularity) of $\overline{X}$ was established in \cite[Proposition 3.8]{CKSV25}. As a consequence it followed that $\overline{X}$ can be refined to be a strongly Feller process starting from every point in $\overline{D}$. The following sharp two-sided estimates of the transition density (i.e., heat kernel) $\overline p$ of  
$\overline{X}$ constitute the main result, Theorem 1.2, of \cite{CKSV25}. 
For simplicity, we state it below for Lipschitz 
open sets, whose boundaries have Assouad dimension  $d-1$, although \cite{CKSV25} allows for more general open sets. 
Note that if $D$ is a bounded Lipschitz open set, or a half-space-like 
Lipschitz open set, or an exterior Lipschitz open set, then ${\rm dim}_{\rm A}(\partial D)=d-1$.

\begin{thm}\label{t:CKSV25}
Suppose $D\subset \R^d$, $d\ge 1$, is a Lipschitz open set.  Assume that ${\rm dim}_{\rm A}(\partial D)= d-1$ if $D$ is unbounded. If \hyperlink{A}{{\bf (A)}}  holds, then 	the process $\overline X$ has a jointly continuous  	transition density  $\overline p(t,x,y)$  defined on $(0,\infty) \times \overline D \times \overline D$. 
Moreover, for any $T \in (0, \infty)$, there exists a constant 
	$C=C(T)\ge 1$ such that for all $(t,x,y) \in (0, T] \times  D\times D$,
	\begin{align*}
	&C^{-1} \left[t^{-d/\alpha} \wedge  \left( \frac{ t }{|x-y|^{d+\alpha}}\Phi\left( \frac{(|x-y| \wedge A_0)^2}{((\delta_D(x) \vee t^{1/\alpha}) \wedge A_0)((\delta_D(y) \vee t^{1/\alpha}) \wedge A_0) }\right) \right)\right]\nn\\
&\le \overline p(t,x,y) \le 
	C \left[t^{-d/\alpha} \wedge  \left( \frac{ t }{|x-y|^{d+\alpha}}\Phi\left( \frac{(|x-y| \wedge A_0)^2}{((\delta_D(x) \vee t^{1/\alpha}) \wedge A_0)((\delta_D(y) \vee t^{1/\alpha}) \wedge A_0) }\right) \right)\right].	\end{align*} 
	\end{thm}
	
\medskip
The main goal of this paper is to establish sharp two-sided heat kernel estimates for processes obtained from $\overline X$ by killing upon approaching the boundary or by a critical potential. 
The framework for this setting was established in the companion paper \cite{CKSV26a}, building on some previous results, primarily \cite{KSV25} where the case of the half-space was studied.

We introduce now the additional assumptions used in \cite{CKSV26a} and refer the reader to 
\cite{CKSV26a}
for motivation and further explanations.

The first assumption concerns the smoothness of the open set. 
\begin{defn}\label{df-c}
	Let $\eta \in (0, 1]$. We say that $D\subset \R^d$, $d\ge 2$, is a $C^{1, \eta}$  open set with characteristics  
	$(\wh R, \Lambda)$, if for any $Q \in \partial D$, there exist a $C^{1, \eta}$ function $\Psi=\Psi^Q:\R^{d-1} \to \R$ with 
	\begin{align*}
		\Psi(\wt 0)= |\nabla \Psi(\wt 0)|=0 \quad \text{and} \quad |\nabla \Psi(\wt y)-\nabla\Psi(\wt z)| 
		\le \Lambda |\wt y-\wt z|^{\eta} \; \text{ for all } \,  \wt y, \wt z \in \R^{d-1},
	\end{align*}
	 and an orthonormal coordinate system {\rm CS}$_{Q}$ with origin at $Q$ such that
	\begin{align}\label{e:local-coordinate_2}
		 B_D(Q, \wh R) = \left\{ y= (\wt y, y_d) \in B(0, \wh R) \text{ in {\rm CS}$_{Q}$} : y_d>\Psi(\wt y) \right\}.
	\end{align}
In the case $d=1$, a $C^{1,\eta}$ open set $D$ with characteristics $(\wh R, \Lambda)$ is a countable union of disjoint open intervals, each of length at least $\wh R$ and the distance between any two distinct intervals being at least $\wh R$. 
\end{defn}

\emph{Throughout this paper, we assume that $D$ is 
a $C^{1,\eta}$ open set with characteristics $(\wh R,\Lambda)$ such that 
${\rm dim}_{\rm A}(\partial D)=d-1$. Without loss of generality, we also assume $\wh R \le 1 \wedge (1/(4\Lambda))^{1/\eta}$. }

Let $\sF^0$ be the closure of  $\mathrm{Lip}_c(D)$ in $L^2(D)$ under $\sE^0_1$. Then $(\sE^0,\sF^0)$ is a regular Dirichlet form. Let $X^0=(X^0_t,t \ge 0; \P_x,x \in D)$ be the Hunt process associated with $(\sE^0,\sF^0)$; see \cite{CKSV26a}. 
When $\alpha\le 1$, the process $\overline X$, started in $D$, never hits $\partial D$ and hence coincides with $X^0$ (see \cite[Proposition 4.9]{CKSV26a}).

The next two conditions,  \hyperlink{B1}{{\bf (B1)}} and \hyperlink{B2}{{\bf (B2)}}, concern  the function $\sB(x,y)$. 
Let $\bH_{-1}:=\{(\wt y, y_d) \in \R^d: y_d>-1\}$ and  define a thin paraboloid-like set $E^Q_\nu(r)$  
as follows: For $Q\in \partial D$, $\nu\in (0,1]$ and $r\in (0, \wh R/4]$, we define
\begin{align}\label{e:EQ}
	E^Q_\nu(r):= \left\{ y=(\wt y, y_d) \text{ in CS$_{Q}$}:  |\wt y|<r/4, \,    4r^{-\nu}|\wt y|^{1+\nu} <y_d < r/2 \right\}.
\end{align}
In the case $d=1$, we interpret $E^Q_\nu(r)$ as the interval $(0,r/2)$.
\medskip

\noindent\hypertarget{B1}{{\bf (B1)}} If $\alpha \ge 1$, then there exist constants $\theta_0>\alpha-1$ and $C_2>0$ such that
\begin{align*}
	|\sB(x,x)-\sB(x,y)| \le C_2 \bigg(\frac{|x-y|}{\delta_D(x) \wedge \delta_D(y) \wedge  \wh R} \bigg)^{\theta_0} 
	\quad \text{for all }  x,y \in D.
\end{align*}
\medskip

\noindent 
\hypertarget{B2}{{\bf (B2)}}
There  exist  constants 
$ \nu_0 \in  (0,\eta)$,   
$\theta_1>0$, $\theta_2>0$,  $C_3>0$ and a non-negative Borel function $\F_0$ on 
$\bH_{-1}$   such that  for any $Q \in \partial D$ and    $x,y \in E^Q_{\nu_0}(\wh R/8)$ with $x=(\wt x,x_d)$ in CS$_Q$,
\begin{align*}
	& \big|\sB(x,y)- \sB(x,x)\F_0((y-x)/x_d) \big| + \big|\sB(x,y)- \sB(y,y)\F_0((y-x)/x_d) \big|\nn \\	
	&\le C_3\left(\frac{|x-y|\vee \delta_D(x)\vee \delta_D(y)}{|x-y|\wedge \delta_D(x)\wedge \delta_D(y)}\right)^{\theta_1} \, 
	(|x-y|\vee\delta_D(x)\vee \delta_D(y))^{\theta_2}.
\end{align*}

\medskip

The interested reader is referred to \cite{CKSV26a} for detailed explanations of \hyperlink{B1}{{\bf (B1)}} and \hyperlink{B2}{{\bf (B2)}}. 
We emphasize that we do not use these two assumptions explicitly in this paper. Instead, we rely solely on their consequences—most notably Theorem \ref{t:Dynkin-improve} (see Section \ref{s:prelim}).

The last condition concerns  the killing potential $\kappa$ and its relation to the function $\F_0$ appearing in \hyperlink{B2}{{\bf (B2)}}.  
Define a  function $\F$ on $\bH_{-1}$ by 
\begin{align}\label{e:def-F0-transform}
	\F(y) =\frac{\F_0(y) + \F_0(-y/(1+y_d))}{2}, \quad \;\;y  = (\wt y,y_d) \in \bH_{-1}.
\end{align}
Recall that $\ub $ is the upper Matuszewska index of the function $\Phi$. 
For  $q \in [(\alpha-1)_+, \alpha- \overline{\beta})$, we define the constant $C(\alpha,q,\F)$ associated with the function $\F$ in \eqref{e:def-F0-transform} by
\begin{align}\label{e:def-C(alpha,q,F)}
	C(\alpha,q,\F) =\begin{cases}	 
		\displaystyle \int_{\R^{d-1}}\frac{1}{(|\wt u|^2+1)^{(d+\alpha)/2}}\int_0^1  
		\frac{(s^q -1)(1-s^{\alpha-1-q})}{(1-s)^{1+\alpha}} \F\big(((s-1)\wt u, s-1)\big) ds \, d \wt u,
		&  d\ge 2, \\[9pt]
		\displaystyle \int_0^1 \frac{(s^q-1)(1-s^{\alpha-q-1})}{(1-s)^{1+\alpha}} \F(s-1 )ds,  & d=1.
	\end{cases}
\end{align}

\medskip
\noindent  \hypertarget{K}{{\bf (K)}} 
There exist constants  
$\theta_3 >0$, $C_4>0$ and 
\begin{align}\label{e:K-upper-bound}	
0\le 	\kappa_0< \lim_{q \uparrow \alpha- \overline{\beta}} C(\alpha,q,\F)
\end{align}such that for all $x \in D$,
	\begin{align*}	\begin{cases}			
	|\kappa(x) - \kappa_0 \sB(x,x) \delta_D(x)^{-\alpha}| \le C_4\delta_D(x)^{-\alpha+\theta_3} &\text{ if }  \delta_D(x) < 1,\\[4pt]			\kappa(x) \le C_4 &\text{ if } \delta_D(x) \ge 1.		
	\end{cases}	
	\end{align*} 		
When $\alpha \le 1$, we further assume that $\kappa_0>0$.

\medskip

 Note that condition \hyperlink{K}{{\bf (K)}} implies that when $\kappa_0 > 0$, the killing potential is  critical: 
 near the boundary, its primary contribution comes from the term $\kappa_0 \sB(x,x) \delta_D(x)^{-\alpha}$. Requiring $\kappa_0 > 0$ when $\alpha \le 1$ ensures that the killing potential $\kappa$ 
 is non-trivial in this case.
 
 It is shown in \cite[Lemma 3.3]{CKSV26a} that	$q \mapsto C(\alpha,q, \F)$ 
 is a well-defined strictly increasing continuous function 
 on $[(\alpha-1)_+, \alpha- \overline{\beta})$ with  $C(\alpha,(\alpha-1)_+,\F)=0$.
  Therefore, under \eqref{e:K-upper-bound},  there exists a unique positive constant 
 $p\in [(\alpha-1)_+, \alpha - \overline{\beta})$ such that 
 \begin{align}\label{e:C(alpha,p,F)}
 	\kappa_0=C(\alpha,p, \F).
 \end{align}
 We note that \eqref{e:K-upper-bound} is not particularly restrictive. For instance, 
 when $\Phi(r)=(r\vee 1)^{{\overline \beta}}(\log (e+r))^a$ for some $a\ge 0$, 
 we have $\lim_{q\uparrow \alpha-\overline\beta}C(\alpha, q, \F)=\infty$, so that \eqref{e:K-upper-bound} holds trivially.

\bigskip

Consider a symmetric form $(\sE^\kappa, \sF^\kappa)$ defined by
	\begin{align*}
		\sE^\kappa(u,v)&=\sE^0(u,v) + 
		\int_{D} u(x)v(x)\kappa(x)dx, \\
		\sF^\kappa&= 	\wt \sF^0 \cap L^2(D, 
			\kappa(x) dx),\nn
	\end{align*}
	where $\wt \sF^0$ is the family of all $\sE^0_1$-quasi-continuous functions in $\sF^0$. Then $(\sE^\kappa,\sF^\kappa)$ is a regular Dirichlet form on $L^2(D)$ with Lip$_c(D)$ as a special standard core,  see \cite[Theorems 6.1.1 and 6.1.2]{FOT}.
Let  $X^\kappa=(X^\kappa_t,t \ge 0; \P_x,x \in D)$ be the Hunt process associated with 
$(\sE^\kappa, \sF^\kappa)$; see \cite{CKSV26a}.  We denote by $\zeta^\kappa$ the lifetime of $X^\kappa$, and define $X^\kappa_t=\partial$ for $t \ge \zeta^\kappa$, where $\partial$ is a cemetery point added to the state space $D$. 
Note that $\kappa$ can be identically zero when $\alpha>1$.	
If $\alpha \le 1$ and $\kappa\equiv 0$, then $X^0= X^\kappa=\overline{X}$ so it is conservative.
Note that when $\alpha \in (1, 2)$, $X^\kappa$ includes $X^0$ as a special case.

The main results of the current paper are the following sharp two-sided estimates of  transition 
density $p^{\kappa}(t,x,y)$
of $X^\kappa$ for small times. 
Note that since $X^{\kappa}$ is a subprocess of $\overline{X}$, it trivially holds that
\begin{align}\label{e:pk-le-pb}
p^{\kappa}(t,x, y)\le \overline{p}(t,x,y)\, .
\end{align}

 Define 
\begin{align}\label{e:def-deltaDt}	\delta_D(x,t)&:=\delta_D(x) \vee t^{1/\alpha}, \quad x\in D, \, t>0.\end{align}

In the following two theorems, we let $p\in [(\alpha-1)_+, \alpha - \overline{\beta})$ be the unique positive constant satisfying \eqref{e:C(alpha,p,F)}.
\begin{thm}\label{t:HKE}
 Suppose $D\subset \R^d$, $d\ge 1$, is a $C^{1,\eta}$ open set. Assume that ${\rm dim}_{\rm A}(\partial D)= d-1$ if $D$ is unbounded. If  conditions \hyperlink{A}{{\bf (A)}}, \hyperlink{B1}{{\bf (B1)}}, \hyperlink{B2}{{\bf (B2)}} and  \hyperlink{K}{{\bf (K)}} hold, then 
the process $X^\kappa$ has jointly continuous transition densities  $p^\kappa(t,x,y)$ defined on $(0,\infty) \times D \times D$.
Moreover, for any  $T>0$, there exist comparison constants depending on $T$ such that the following estimates hold: for all $t\in (0,T]$ and $x,y \in D$,
	\begin{align*}
	&	p^\kappa(t,x,y) \asymp \P_x(\zeta^\kappa>t)   \P_y(\zeta^\kappa>t)\, \overline p(t,x,y)\\
&\asymp \left(1 \wedge \frac{\delta_D(x)}{t^{1/\alpha}}
		\right)^p
		\left(1 \wedge \frac{\delta_D(y)}{t^{1/\alpha}}
		\right)^p
	\left[t^{-d/\alpha} \wedge  \left( \frac{ t }{|x-y|^{d+\alpha}}\Phi\left( \frac{(|x-y| \wedge A_0)^2}{
		(\delta_D(x,t) \wedge A_0)(\delta_D(y,t) \wedge A_0) }\right) \right)\right].
\end{align*}	
\end{thm}

\medskip
Additionally, in the case when $D$ is bounded, we also obtain large-time estimates.  
Note that if $D$ is a bounded  $C^{1,\eta}$ open set,  then ${\rm dim}_{\rm A}(\partial D)= d-1$.

\begin{thm}\label{t:main2}
 Suppose $D\subset \R^d$, $d\ge 1$, is a bounded $C^{1,\eta}$ open set, and that conditions \hyperlink{A}{{\bf (A)}},  \hyperlink{B1}{{\bf (B1)}}, \hyperlink{B2}{{\bf (B2)}} and \hyperlink{K}{{\bf (K)}} hold. Then for any $T>0$, there exist comparison constants depending on $T$ such that the following estimate holds:
 for all $t\ge T$ and $x,y\in D$,
  \begin{align}\label{t:e:Lhke}
 	p^\kappa(t,x,y)
 	& \asymp   e^{-\lambda_1 t}\, \delta_D (x)^{p}\, \delta_D (y)^{p} \asymp  e^{-\lambda_1 t}\, \P_x(\zeta^\kappa>1)  \P_y(\zeta^\kappa>1),
\end{align}
where $-\lambda_1<0$ is the largest eigenvalue of the infinitesimal generator of $X^\kappa$.
\end{thm}

By combining Theorems \ref{t:HKE} and \ref{t:main2}, we establish in Theorem \ref{t:gfe} sharp two-sided estimates of the Green function of $X^{\kappa}$ in a bounded set $D$.

\medskip
\textbf{Organization of the paper.} 
 Section \ref{s:examples} is devoted to examples illustrating the scope and applicability of our results. Readers primarily interested in the general theory may skip this section without affecting their understanding of the subsequent sections.

In Section \ref{s:prelim} we list several results obtained mainly in \cite{CKSV25, CKSV26a}. Among those we single out Theorem \ref{t:Dynkin-improve} that provides the exact decay rate.

Sections \ref{s:uHKE} and \ref{s:lHKE} form the core of the paper. Here we establish the upper and lower bounds on 
$p^{\kappa}(t,x,y)$ for short time. The upper bound is more difficult and requires some novel ideas. This is mainly due to the fact that when $ \overline{\beta} >\alpha/2$,  both $p^\kappa(t,x,y)$ and $\overline p(t,x,y)$ may not be uniformly bounded on $[0, T] \times D \cap\{|x-y|> \delta\}$ for any $\delta>0$. Thus, the well-known general upper bounds like those in 
\cite[Lemma 4.2]{CKS14}, \cite[Lemma 2]{BGR10} and their variations
that were used frequently before, are not sufficient in our case, and are replaced by more delicate bounds obtained in Lemma \ref{l:general-upper-2}.
Another novelty is Proposition \ref{p:exitu} providing negative moment estimates of exit times. The short time upper bound  
 involving only one boundary term is given in Proposition \ref{p:blow-up-step2}. The rather long proof relies on already mentioned Lemma \ref{l:general-upper-2} and Proposition \ref{p:exitu}, as well as several estimates of integrals involving heat kernels, blow-up weight function, and negative powers of distance to the boundary. Combined with the semigroup property and Theorem \ref{t:CKSV25}, the estimate in Proposition \ref{p:blow-up-step2} implies the full short time upper bound given in Proposition \ref{p:HKE-generalu}.

The lower bound for short time is proved in Section \ref{s:lHKE}. 
The main ingredients and ideas used to obtain the short-time lower bound are similar to those in \cite[Section 5]{CKSV25b}, where the corresponding lower bound was established for purely discontinuous Markov processes with jump kernels decaying at the boundary of a bounded $C^{1,1}$ open set. However, to cover unbounded open sets as well, we establish in Proposition \ref{p:ndl2} 
a near-diagonal heat kernel lower bound
that holds uniformly on balls of arbitrary radius away from the boundary.
Combining this interior lower bound with the interior estimates for the mean exit time established in
Lemma \ref{l:E2-lower}, we obtain the sharp lower bound for the mean exit time near the boundary in Proposition \ref{p:met-bound-new}.
Using the results above and the basic properties of lifted points obtained in \cite{CKSV25b}, we obtain short time heat kernel lower bound in Proposition \ref{p:lHKE-small}. Lemmas \ref{l:lift-probability}–\ref{l:lift} and Proposition \ref{p:lHKE-small} are the counterparts, in the present setting, 
of
 Lemmas 5.4–5.6 of \cite{CKSV25b}.

In Section \ref{s:main-proofs}, we first prove Proposition \ref{p:life} which gives sharp two-sided estimates on the survival probability.  
Then combining Propositions \ref{p:HKE-generalu},  \ref{p:lHKE-small} and monotonicity of the heat kernels, we obtain Theorem \ref{t:HKE-general}, a slight generalization of Theorem \ref{t:HKE}. The latter result then follows from Theorem \ref{t:HKE-general} and  Proposition \ref{p:life}.  
The proof of Theorem \ref{t:main2} is quite standard and only a brief outline is provided. These two theorems are then used to prove the sharp two-sided estimates of the Green function for bounded $D$ in Theorem \ref{t:gfe}.

\textbf{Notation.} We end this section with a few words on notation. Throughout this 
work, we use ``$:=$" to denote a definition, which is read as ``is defined to be''. We use the notation
$a\wedge b:=\min \{a, b\}$ and  $a\vee b:=\max\{a, b\}$. $\N$ denotes the set of natural numbers.  
The notation $C=C(a,b,\ldots)$ indicates that the constant $C$
depends on $a, b, \ldots$. 
The dependence on  $d, \alpha$, the localization characteristics $\wh R$ and $\Lambda$ and the constants in conditions \hyperlink{A}{{\bf (A)}}, \hyperlink{K}{{\bf (K)}}, \hyperlink{B1}{{\bf (B1)}} and \hyperlink{B2}{{\bf (B2)}}  may not be mentioned explicitly. 
Constants $\kappa_0$, $\theta_3$ (from \hyperlink{K}{{\bf (K)}}), $\theta_0$ (from \hyperlink{B1}{{\bf (B1)}}), $\nu_0, \theta_1, \theta_2$ (from \hyperlink{B2}{{\bf (B2)}}), $\epsilon_1$ (from Theorem \ref{t:Dynkin-improve}), and $\wh r, \eta_0$ (from \eqref{e:lifting-property}) are fixed throughout the paper. 
Upper case letters $C_i$, $i \in \N$,  with subscripts denote  strictly positive constants whose values appear in definitions and statements of the results, and are fixed throughout this work.
The values of constants $c_i$, $i \in \N_0$, are fixed in each proof, and the labeling of these constants starts anew in each proof. 
We denote $x\in \R^d$ as $x=(\wt{x}, x_d)$ with $\wt{x}\in \R^{d-1}$.
We  use $m_d$ to denote the Lebesgue measure on $\R^d$. 
For a Borel subset $A\subset \R^d$, $\delta_A(x)$ denotes the Euclidean distance between $x$ and $\partial A$. 
For a  subset $A \subset \R^d$,  we define
\begin{align*}	B_{A}(x,r):=A \cap B(x,r), \quad x \in \R^d, \, r>0.
\end{align*}
    We adopt the convention $c/0=\infty$ for $c>0$.


\section{Examples}\label{s:examples}

This section 
explores three examples that fit our framework: 
  elliptic  problems related to  the nonlocal Neumann condition introduced in \cite{DRV17},  path-censored $\alpha$-stable process, and path-censored processes with critical killings in the open half-space discussed in \cite{KSV25}.

Throughout this  
section we assume that $d\ge 1$ and $\alpha \in (0,2)$. Recall that the fractional Laplacian $\Delta^{\alpha/2}:=-(-\Delta)^{\alpha/2}$ on $\mathbb{R}^d$ is defined as
\begin{align*}
	\Delta^{\alpha/2} u(x)&:= p.v. \int_{\R^d} (u(y)-u(x)) J^\alpha(x,y)dy \\
	&:= \lim_{\eps \downarrow 0} \int_{\{y\in \R^d: |x-y|>\eps\}} (u(y)-u(x)) J^\alpha(x,y)dy,
\end{align*}
where 
\begin{align*}	
	J^\alpha(x,y):=A_{d,-\alpha} |x-y|^{-d-\alpha}, \qquad A_{d,-\alpha}:=\frac{2^\alpha \Gamma((d+\alpha)/2)}{\pi^{d/2} |\Gamma(-\alpha/2)|}.\end{align*}

\begin{defn}
	\rm 	We say that an open set $D\subset \R^d$ is \textit{half-space-like} if there are constants $a<b$ such that $\{(\wt x,x_d)\in \R^d: x_d>b\} \subset D \subset \{(\wt x,x_d)\in \R^d: x_d>a\}$.
\end{defn}

\subsection{Elliptic problems related to nonlocal Neumann condition}
This example is a continuation of \cite[Subsection 6.1]{CKSV25} and \cite[Subsection 9.1]{CKSV26a}. Let 
 $D\subset \R^d$ be a $C^{1,\eta}$ open set  that is either bounded,  or half-space-like, or the complement of a bounded set. Note that  ${\rm dim}_{\rm A}(\partial D)=d-1$. 

In \cite{DRV17}, the nonlocal Neumann problem was introduced as 
\begin{align}\label{e:nonlocal-Neumann}
	\begin{split} 
		\Delta^{\alpha/2} u=f \quad  &\text{ in $D$,}\\
		\sN^{\alpha/2}_D u = 0  \quad  &\text{ in $\R^d \setminus D$,} \end{split} 
\end{align}
employing the nonlocal normal derivative $\mathcal{N}^{\alpha/2}_D u(x) = \int_D (u(x)-u(y)) J^\alpha(x,y) dy$. Following \cite{Ab20}, we can rewrite \eqref{e:nonlocal-Neumann} as a regional nonlocal problem in $D$ with 
 the operator 
 $$\Delta^{\alpha/2,N}_D u(x):= p.v. \int_{D} (u(y)-u(x)) K^\alpha_D(x,y)dy,$$  
where the modified kernel $K^\alpha_D(x,y)$ is defined by
\begin{align*}
	K^\alpha_D(x,y):= J^\alpha(x,y) + \int_{D^c} \frac{J^\alpha(x,z)J^\alpha(z,y)}{ \int_D J^\alpha(z,w)dw } \, dz, \quad x,y \in D.
\end{align*}
We refer to \cite{Von21} for further discussion on the probabilistic aspects of this operator.

Let $\sB_D^\alpha(x,y)$ denote the  function associated with $K_D^\alpha(x,y)$, defined by
\begin{align*}
	\sB^\alpha_D(x,y):=\begin{cases}
		K_D^\alpha(x,y)|x-y|^{d+\alpha} &\mbox{ if $x\neq y$},\\[5pt]
		A_{d,-\alpha}&\mbox{ if $x=y$}.
	\end{cases}
\end{align*}
By 
\cite[Proposition 6.3]{CKSV25},
$\sB^\alpha_D(x,y)$ satisfies condition \hyperlink{A}{\bf (A)} with $\Phi(r) = \log(e+r)$ and $A_0=\diam(D^c)$. 
Moreover,
$\sB^\alpha_D(x,y)$ also satisfies \hyperlink{B1}{\bf (B1)} 
by \cite[Section 7.1]{KSV23}
and satisfies \hyperlink{B2}{\bf (B2)} by \cite[Proposition 9.2]{CKSV26a}.
Thus, applying Theorems \ref{t:HKE} and \ref{t:gfe}, we obtain the following explicit examples.  

\begin{example}\label{ex:1}
	\rm {\bf (Nonlocal Neumann problem with boundary Hardy potential)} Let   $D\subset \R^d$ be a $C^{1,\eta}$ open set  that is either bounded,  or half-space-like, or the complement of a bounded set. Let $\kappa \in [0,\infty)$. When $\alpha\le 1$, we assume that $\kappa>0$. 
	Consider the  following nonlocal Neumann problem on $D$ with a boundary Hardy potential:
	\begin{align*}
		\begin{cases} 
			{\Delta^{\alpha/2} u(x)=\kappa \delta_D(x)^{-\alpha}u(x) + f(x)}  & \text{ in $D$,}\\[1pt]
			u = 0 &  \text{ on $\partial D$,}\\[1pt]
			{\sN^{\alpha/2}_D u =0}    &\text{ in $\R^d \setminus 
				\overline D$,} 
		\end{cases}    
	\end{align*}
which is associated with the following   variational formulation:
	\begin{align*}
		\sE^{\alpha, N,\kappa}(u,v)&:= \frac12 \int_{D\times D} (u(x)-u(y))(v(x)-v(y))K_D^\alpha(x,y)dxdy   + \kappa \int_{D} u(x) v(x) \delta_D(x)^{-\alpha}dx,\\
		\sF^{\alpha,N,\kappa}&:=\text{$(	
			\sE^{\alpha, N,\kappa})_1
			$-closure of $C_c^2(D)$ in $L^2(D)$}.
	\end{align*}
	Since $(\sE^{\alpha, N,\kappa}, \sF^{\alpha,N,\kappa})$ falls within our framework, we obtain from
	Theorems \ref{t:HKE} and \ref{t:gfe} that
	\begin{thm} \label{t:kernelestimate-Neumann}
		 Suppose that  $D$ is $C^{1,\eta}$, and is either bounded or half-space-like or the complement of a bounded set.  Then  $(\sE^{\alpha, N,\kappa}, \sF^{\alpha,N,\kappa})$ has a jointly continuous heat kernel $p(t,x,y)$ on $(0,\infty) \times D\times D$.  Moreover, there exists a unique $p_\kappa \in [(\alpha-1)_+, \alpha) \cap  (0,\alpha)$  such that 
		  for every $T>0$,  the following estimates hold for $(t,x,y) \in (0,T] \times D\times D$: 
			\begin{align*}		&p(t,x,y)   \asymp \left(1 \wedge \frac{\delta_D(x)}{t^{1/\alpha}}		\right)^{p_\kappa}		\left(1 \wedge \frac{\delta_D(y)}{t^{1/\alpha}}		\right)^{p_\kappa} 		\\		&\;\; \times \left[t^{-d/\alpha} \wedge		\left(\frac{ t }{|x-y|^{d+\alpha}}\log\left( e+\frac{(|x-y| \wedge \diam(D^c))^2}{((\delta_D(x) \vee t^{1/\alpha}) \wedge \diam(D^c))((\delta_D(y) \vee t^{1/\alpha}) \wedge \diam(D^c)) }\right)\right)\right].	\end{align*} 	Furthermore, if $D$ is bounded, then \eqref{e:FGE} holds true with $p=p_\kappa$.\end{thm}
	
	We note that when $\alpha>1$ and $\kappa=0$,  the constant $p_\kappa$ in Theorem \ref{t:kernelestimate-Neumann} is given by  $p_\kappa=\alpha-1$.
\end{example}

  Define
\begin{align*}
	\kappa_D(x):=\int_{D^c} J^\alpha(x,y)dy, \quad x\in D,\\
		\wh \kappa_D(z):= \int_DJ^\alpha(z,y)dy, \quad z \in \overline D^c.
\end{align*}

\begin{example}\label{ex:2}
	\rm {\bf (Nonlocal Robin problem)} Let $\theta \in (0,1)$.  Consider the elliptic problem with a nonlocal Robin boundary condition (see \cite{Ab20}):
	\begin{align}\label{e:nonlocal-Robin}
		\begin{cases} 
			{\Delta^{\alpha/2} u=f}  & \text{ in $D$,}\\[1pt]
			{\theta \wh \kappa_D(z) u(z)+  (1-\theta)\sN^{\alpha/2}_D u(z) =0}  & \text{ in $\R^d \setminus 
				\overline D$.} 
		\end{cases}    
	\end{align}
	In view of \cite[Theorem 1.1]{Ab20},  the variational formulation of \eqref{e:nonlocal-Robin} is given by
	\begin{align*}
		\sE^{\alpha, R,\theta}(u,v)&:= \frac12 \int_{D\times D} (u(x)-u(y))(v(x)-v(y))K_D^{\alpha,\theta}(x,y)dxdy   + \theta \int_{D} u(x) v(x) \kappa_D(x)dx,\\
		\sF^{\alpha,R,\theta}&:=\text{$(	\sE^{\alpha, R,\theta})_1$-closure of $C_c^2(D)$ in $L^2(D)$}.
	\end{align*}
	Note that	the  function 
	$\sB_D^{\alpha,\theta}(x,y)$ associated with $K_D^{\alpha,\theta}(x,y)$ is  
	\begin{align*}
		\sB^{\alpha,\theta}_D(x,y)=\theta A_{d,-\alpha} + (1-\theta) \sB^\alpha_D(x,y), \quad x,y \in D.
	\end{align*} 
	
	Since 
	$\sB^\alpha_D(x,y)$ satisfies condition \hyperlink{A}{\bf (A)} with $\Phi(r) = \log(e+r)$ and $A_0=\diam(D^c)$, $\sB^{\alpha,\theta}_D(x,y)$ also satisfies condition \hyperlink{A}{\bf (A)} with the same $\Phi$ and $A_0$. Further, by \cite[Proposition 9.2 and Lemma 10.6]{CKSV26a},  $\sB_D^{\alpha,\theta}(x,y)$ satisfies  \hyperlink{B1}{{\bf (B1)}} and \hyperlink{B2}{{\bf (B2)}}, and the killing potential $\theta\kappa_D(x)$ satisfies  \hyperlink{K}{{\bf (K)}}. Hence, we deduce that
		\begin{thm} \label{t:kernelestimate-Neumann-2}
		Suppose that  $D$ is $C^{1,\eta}$, and is either bounded or half-space-like or the complement of a bounded set.  Then  $(\sE^{\alpha, R,\theta}, \sF^{\alpha, R,\theta})$ has a jointly continuous heat kernel $p(t,x,y)$ on $(0,\infty) \times D\times D$.  Moreover, there exists a unique constant $p(\theta) \in [(\alpha-1)_+, \alpha) \cap  (0,\alpha)$  such that 
		for every $T>0$,  the following estimates hold for $(t,x,y) \in (0,T] \times D\times D$: 
		\begin{align*}		&p(t,x,y)   \asymp \left(1 \wedge \frac{\delta_D(x)}{t^{1/\alpha}}		\right)^{p(\theta)}		\left(1 \wedge \frac{\delta_D(y)}{t^{1/\alpha}}		\right)^{p(\theta)} 		\\		&\;\; \times \left[t^{-d/\alpha} \wedge		\left(\frac{ t }{|x-y|^{d+\alpha}}\log\left( e+\frac{(|x-y| \wedge \diam(D^c))^2}{((\delta_D(x) \vee t^{1/\alpha}) \wedge \diam(D^c))((\delta_D(y) \vee t^{1/\alpha}) \wedge \diam(D^c)) }\right)\right)\right].	\end{align*} 	Furthermore, if $D$ is bounded, then \eqref{e:FGE} holds true with $p=p(\theta)$.\end{thm}
	
	Note that as $\theta \to 0$ and $\theta \to 1$, respectively, the corresponding limit problems of the Robin problem \eqref{e:nonlocal-Robin} are the problem with $\kappa=0$ considered in Example \ref{ex:1} and the nonlocal Dirichlet problem. Hence, 
	since $\sB_D^{\alpha,\theta}(x,y)$ is decreasing in $\theta$ and the killing potential $\theta \kappa_D(x)$ is increasing in $\theta$, one can verify that $p(\theta)$ is increasing in $\theta$ and $p(\theta) \in [(\alpha-1)_+, \alpha/2) \cap  (0,\alpha/2)$. Moreover, since the boundary decay rate for the Dirichlet problem is known to be $\delta_D(x)^{\alpha/2}$, we have
  $\lim_{\theta\to 0} p(\theta) = (\alpha-1)_+$ and $\lim_{\theta \to 1} p(\theta) = \alpha/2$. 
\end{example}

\subsection{Path-censored $\alpha$-stable process in open sets}

This example is a continuation of \cite[Subsection 6.2]{CKSV25} and \cite[Subsection 9.2]{CKSV26a}. Suppose that  $D\subset \R^d$ is a $C^{1,\eta}$ open set  that is either bounded,  or half-space-like, or the complement of a bounded set. 
Let $Y=(Y_t)_{t\ge 0}$ be an isotropic $\alpha$-stable process in $\R^d$ with generator $\Delta^{\alpha/2}$. Define
$$\sA^{\overline D}_t:=\int_0^t \1_{\{Y_s \in \overline D\}}ds,\quad t\ge 0,$$ 
and let $\tau^{\overline D}_t:=\inf\{s>0:\sA^{\overline D}_s>t\}$ be the right-continuous inverse of $\sA^{\overline D}$. The \textit{trace process} $X$ of $Y$ on $\overline D$ is defined by the time change $X_t:=Y_{\tau^{\overline D}_t}$. We refer to \cite[Chapter 5]{CF12} for basic results on trace processes.
The subprocess $X^D$ 
obtained by killing $X$ upon leaving $D$, or equivalently, hitting $\partial D$,  is referred to as the \textit{path-censored $\alpha$-stable process} on $D$; see \cite{KPW14}.   For $d=1$, path-censored $\alpha$-stable processes were studied in  \cite{KPW14}.  
When $\alpha\le 1$, by \cite[Proposition 4.9]{CKSV26a},  $\partial D$ is polar for $X$, so that $X^D$ coincides with $X$, whose heat kernel estimates were studied in \cite[Theorem 1.2 and Subsection 6.2]{CKSV25}. 
 In this subsection we consider the case $d\ge 2$ and $\alpha \in (1,2)$.

By \cite[Proposition 6.4]{CKSV25}, the path-censored $\alpha$-stable process  satisfies condition \hyperlink{A}{\bf (A)} 
with $\Phi(r) =1+ r^{\alpha/2}$ and $A_0=\diam(D^c)$.  It also satisfies \hyperlink{B1}{\bf (B1)}  by \cite[Section 7.1]{KSV23} and  \hyperlink{B2}{\bf (B2)} by \cite[Proposition 9.5]{CKSV26a}. When $D$  is either  half-space-like or the complement of a bounded set, the Dirichlet form of $X^D$ is purely jump, with no killing part.  When $D$ is bounded, the Dirichlet form of $X^D$ 
 has a non-trivial  killing measure with density $\kappa(x)$  comparable to $\delta_D(x)^{-\alpha/2}$; see \cite[Subsection 6.2]{CKSV25}.  Hence, in all cases, \hyperlink{K}{{\bf (K)}} holds with $\kappa_0=0$ and $\theta_3=\alpha/2$. Consequently, applying Theorems \ref{t:HKE} and \ref{t:gfe}, we obtain the following result.

\begin{thm}\label{t:kernelestimate-trace}
	Suppose that $d\ge 2$, $\alpha \in (1,2)$, and that $D$ is $C^{1,\eta}$, and is either bounded or half-space-like or the complement of a bounded set.
	Then the  path-censored $\alpha$-stable process $X^D$ on $D$ has a jointly continuous heat kernel $p(t,x,y)$ on $(0,\infty) \times D\times D$ 
	and, for every $T>0$, it satisfies  the following estimates for $(t,x,y) \in (0, T] \times D\times D$: 
	If $|x-y|\le t^{1/\alpha}$, then
	\begin{align*}
		p(t,x,y) \asymp
		\left(1 \wedge \frac{\delta_D(x)}{t^{1/\alpha}}
		\right)^{\alpha-1}
		\left(1 \wedge \frac{\delta_D(y)}{t^{1/\alpha}}
		\right)^{\alpha-1} t^{-d/\alpha},
	\end{align*}
	and if $|x-y|>t^{1/\alpha}$, then
	\begin{align*}
			p(t,x,y) 
		&\asymp \left(1 \wedge \frac{\delta_D(x)}{t^{1/\alpha}}
		\right)^{\alpha-1}
		\left(1 \wedge \frac{\delta_D(y)}{t^{1/\alpha}}
		\right)^{\alpha-1}\frac{ t }{|x-y|^{d+\alpha}}
		\\
		&\quad \times
		\left[ 1+  \bigg( \frac{(|x-y|\wedge \diam(D^c))^2}{((\delta_D(x)\vee t^{1/\alpha})\wedge \diam(D^c))((\delta_D(y)\vee t^{1/\alpha})\wedge \diam(D^c))}\bigg)^{\alpha/2}\right] .
	\end{align*} 
	Furthermore, if $D$ is bounded, then \eqref{e:FGE} holds true with $p=\alpha-1$.
\end{thm}

\subsection{Path-censored self-similar Markov processes with critical killing in the open half-space}\label{ss:res-proc}
We briefly recall the construction of resurrected processes in the upper half-space 
 as detailed in \cite[Section 2]{KSV25}. 
We will use the resurrected process to obtain a path-censored process $X^\kappa$ with critical killing $\kappa x_d^{-\alpha}$ that fits within our framework.

In this setting, an $\alpha$-stable process starting in 
 the upper half-space $\bH$ is {\it resurrected} whenever it exits to the lower half-space 
$-\bH:=\{x=(\wt{x}, x_d): \wt{x}\in \R^{d-1}, x_d<0\}$. Specifically, if the process lands at $z \in -\bH$, it is immediately relocated to a point $y \in \bH$ according to probability distribution 
$\mu(z,y)dy$. The interaction between the starting point $x$ and the return point $y$ is captured by the resurrection kernel $q^\alpha_\mu(x,y)$:
$$
q^\alpha_\mu(x,y) = \int_{-\bH} J^\alpha(x,z) \mu(z,y) \, dz, \quad x,y \in \bH.
$$
Consequently, the jump kernel for the resurrected process incorporates these additional transitions, taking the form:
\begin{equation}\label{e:jump-kernel-res}
	J^{\alpha}_{\mu}(x,y)=J^\alpha(x,y) +q^\alpha_\mu(x,y), \quad x,y \in \bH.
\end{equation}
 Define 
\begin{align*}
	\sB^\alpha_\mu(x,y):=\begin{cases}
		J^\alpha_\mu(x,y)|x-y|^{d+\alpha} &\mbox{ if $x\neq y$},\\[5pt]
		A_{d,-\alpha}&\mbox{ if $x=y$},
	\end{cases}\qquad x,y \in \bH.
\end{align*}
Note that for certain values of $\alpha \in (0,2)$, the process may be conservative or exhibit boundary killing at $\partial \bH$.

Let $\Psi$ be an almost increasing positive function on $[0, \infty)$ 
satisfying  the following weak scaling condition on $[1,\infty)$: there exist constants $C>0$ and $\gamma\in [0, 1 \wedge \alpha)$ such that 
$$
\frac{\Psi(r)}{\Psi(s)}\leq C \bigg(\frac{r}{s}\bigg)^{\gamma}\quad  \text{for all $1\le s\le r$.}
$$
Assume that the return probability density $\mu(z,y)$ has the form  
$$
 \mu(z,y)=C_\Psi{|z_d|^{\alpha}}\Psi\left(\frac{|y-z|^2}{y_d|z_d|}\right)|y-z|^{-d-\alpha}, \quad z\in -\bH, \, y \in \bH,
$$
where $C_\Psi$ is the normalizing constant. 
It is shown in \cite[Proposition 2.2]{KSV25} that 
the resurrection kernel $q_\mu(x,y)$ 
is symmetric for this kind of $\mu$.
We note in passing that when  
$\Psi(r)=r^{\alpha/2}$, the resurrected process coincides with the trace process 
(killed upon hitting $\partial \bH$ when $\alpha\in (1,2)$),
while when $\Psi(r)=1$, we recover the process associated with nonlocal Neumann problem 
(killed upon hitting $\partial \bH$ when $\alpha\in (1,2)$). 

Define $\Psi_1(u):=1$ for $u\in [0,2]$ and 
$$
\Psi_1(u):=\left(\int_1^2 \frac{\Psi(v)}{v}\, dv \right)^{-1}\int_1^u \frac{\Psi(v)}{v}\, dv\quad \text{for $u \in [2,\infty)$} .
$$
It is shown in \cite[Theorem 2.4]{KSV25} that  
$\sB^\alpha_\mu$ satisfies assumptions \hyperlink{A}{\bf (A)} and \hyperlink{B1}{\bf (B1)} 
with the 
blow-up weight function $\Phi=\Psi_1$ and $A_0=\infty$. Moreover, by
the homogeneity and  the horizontal translation invariance of $\mu$ in \cite[(2.1)-(2.2)]{KSV25}, 
$\sB^\alpha_\mu$  
is   invariant under  dilations and horizontal translations. That is, 
\begin{equation}\label{e:invariant-properties}
	\sB^\alpha_\mu(x,y)  = \sB^\alpha_\mu(ax,ay) =  \sB^\alpha_\mu(x+ (\wt z,0), y+(\wt z,0)) \quad \text{for all $x,y\in \bH$, $a>0$, and $\wt z\in\R^{d-1}$}.
\end{equation}
Define $\F^\alpha_\mu(y):=A_{d,-\alpha}^{-1}\sB^\alpha_\mu(\e_d, y+\e_d)$ for $y \in \bH_{-1}$. By \eqref{e:invariant-properties}, we have 
	$\sB^\alpha_\mu(x,y) = \sB^\alpha_\mu(x,x)\F^\alpha_\mu((y-x)/x_d)$ for all $x,y \in \bH$, and $ \F^\alpha_\mu(y) =  \F^\alpha_\mu( - y/(y_d+1))$ for all $y \in \bH_{-1}$. See the discussion at the beginning of \cite[Section 3]{CKSV26a} for details.  Consequently,  \hyperlink{B2}{\bf (B2)} holds trivially with $\F_0= \F^\alpha_\mu$. 
For $\kappa \ge 0$, we consider the bilinear form 
$$
\EE^{\alpha,\mu,\kappa}(u,v)=\frac{1}{2}\int_{\bH \times \bH} (u(x)-u(y))(v(x)-v(y))J^\alpha_\mu(x,y)\, dx, dy+  \kappa \int_{\bH} u(x)v(x) x_d^{-\alpha}dx,
$$
where the jump kernel $J^\alpha_\mu(x,y)$ is given by \eqref{e:jump-kernel-res}.
Let $\FF^{\alpha,\mu,0}$ be the closure of $C_c^{\infty}(\bH)$ in $L^2(\bH,dx)$ under the norm $(\sE^{\alpha,\mu,0}_1)^{1/2}$, where $\sE^{\alpha,\mu,0}_1:= \sE^{\alpha,\mu,0} + \lVert \cdot \rVert_{L^2(\bH)}^2$.
For $\kappa>0$, define 
$$
\sF^{\alpha,\mu,\kappa}:= 	\wt\sF^{\alpha,\mu,0} \cap L^2(\bH,  x_d^{-\alpha}dx),\nn
$$
where $\wt\sF^{\alpha,\mu,0}$ is the family of all $\sE^{\alpha,\mu,0}_1$-quasi-continuous functions in $\sF^{\alpha,\mu,0}$.  
It is shown in \cite{KSV25} that 
$(
\sE^{\alpha,\mu,\kappa},
\sF^{\alpha,\mu,\kappa})$ is  a regular Dirichlet form on $L^2(\bH, dx)$ and that its associated Hunt process  $X^{\mu,\kappa}=((X_t^{\mu,\kappa})_{t\ge 0}, (\P_x)_{x\in \bH})$ is self-similar and invariant under horizontal translations. 

Let $\ub$ be the upper Matuszewska index of  
$\Psi_1$. For $q\in [(\alpha-1)_+, \alpha-\ub)$, define
\begin{align*}
&	C(\alpha,q,\mu)=\begin{cases}	 
	\displaystyle \int_{\R^{d-1}}\frac{1}{(|\wt u|^2+1)^{(d+\alpha)/2}}\int_0^1  
	\frac{(s^q -1)(1-s^{\alpha-1-q})}{(1-s)^{1+\alpha}} \sB^\alpha_\mu( \e_d, ((s-1)\wt u, s)  ) ds \, d \wt u,
	&  d\ge 2, \\[9pt]
	\displaystyle \int_0^1 \frac{(s^q-1)(1-s^{\alpha-q-1})}{(1-s)^{1+\alpha}} \sB^\alpha_\mu(1,s)ds,  & d=1.
\end{cases}
\end{align*}
Note that $C(\alpha,q,\mu) = A_{d,-\alpha} C(\alpha,q,\F^\alpha_\mu)$, where $C(\alpha,q,\F^\alpha_\mu)$ is the constant associated with $\F^\alpha_\mu$ defined by \eqref{e:def-C(alpha,q,F)}.
By Theorem \ref{t:HKE} and 
the self-similarity of $X^{\mu,\kappa}$, 
we obtain the following heat kernel estimates for all time.  
\begin{thm}\label{t:kernelestimate-res}
Let $d\ge 1$ and $\alpha\in (0,2)$. Suppose that $0\le \kappa < \lim_{q\uparrow \alpha-\ub}C(\alpha,q,\mu)$, and assume further that $\kappa>0$ when $\alpha\le 1$. 
Then the process $X^{\mu,\kappa}$  has a jointly continuous heat kernel $p(t,x,y)$ on $(0,\infty) \times \bH \times \bH$ 
	that satisfies  the following estimates for all $(t,x,y) \in (0, \infty) \times \bH \times \bH$: 
	\begin{align*}
		p(t,x,y) 
		&  \asymp \left(1 \wedge \frac{
			 x_d}{t^{1/\alpha}}
		\right)^{p_\kappa}
		\left(1 \wedge \frac{
			 y_d}{t^{1/\alpha}}
		\right)^{p_\kappa} \left[t^{-d/\alpha} \wedge  \left(\frac{ t }{|x-y|^{d+\alpha}}\Psi_1\left( \frac{|x-y|^2}{(x_d\vee t^{1/\alpha})(y_d\vee t^{1/\alpha})}\right)\right)
		\right],
	\end{align*} 
	where $p_\kappa \in [(\alpha-1)_+, \alpha-\ub)$ is the unique positive constant such that $\kappa = C(\alpha, p_\kappa, \mu)$.
\end{thm}

\section{Preliminary results}\label{s:prelim}

In this section, we collect several preliminary results, mainly from \cite{CKSV25} and \cite{CKSV26a}.

 Recall that the function $\Phi$ has  upper Matuszewska index $\ub$.  By definition, for any $\beta>\ub$,
 there is a constant 
 $C=C(\beta)\ge 1$  such that the following holds:
\begin{align}\label{e:Phi-scaling}
	\frac{\Phi(r)}{\Phi(s)}\le C\bigg(\frac{r}{s}\bigg)^{\beta} 
	\quad \text{for all $0<s\le r$}.
\end{align} 

The following elementary inequality will be used frequently:
\begin{align}
	\frac{(r \vee B)\wedge A}{(s \vee B) \wedge A} \le \frac{r}{s} \quad \text{for all $B \in [0,\infty)$, $A\in (0,\infty]$ and $0<s\le r$}.\label{e:without-A0}
\end{align}

The following lemma is a consequence of  \eqref{e:Phi-scaling} and \eqref{e:without-A0}; see the proof of  \cite[Lemma 2.8]{CKSV25}.
\begin{lemma}\label{l:jump-density-monotonicty}
	There exists $C>0$ such that for all $a, r>0$ and $A \in (0, \infty]$,  
	\begin{align*}
		&	\sup_{R\ge r}	\frac{ 1 }{R^{d+\alpha}}\Phi\left( \frac{(R \wedge A)^2}{a}\right)
		\le 	\frac{ C }{r^{d+\alpha}}\Phi\left( \frac{(r\wedge A)^2}{a}\right).
	\end{align*}
\end{lemma}

Recall that $B_D(x,r)=D\cap B(x,r)$. 
The next  lemma, which is essentially due to \cite[Lemma 2.6]{CKSV25},   provides integral estimates for powers of the boundary distance function. These estimates will be used repeatedly to control the blow-up behavior of the jump kernel near the boundary.

\begin{lemma}\label{l:Aikawa-integrals}
	For any $q >-1$,	there exists $C=C(q)>0$ such that
	\begin{align*}
		\int_{B_D(x,r)} {(\delta_D(y) \wedge A)^{q}}   dy\le Cr^{d}((\delta_D(x) \vee r) \wedge A)^{q}   \quad \text{ for all $A\in 		(0,\infty]$, $x\in  {\overline D}$ and  $r>0$.}
	\end{align*}
\end{lemma}
\begin{proof}
	When $q\in (-1,0]$, the result was established in \cite[Lemma 2.6(i)]{CKSV25} (note that $\gamma$ there is equal to 1 in the current setting).  For the case $q>0$, the result follows from the fact that $\delta_D(y) \le \delta_D(x) +r \le 2(\delta_D(x)\vee r)$ for all $y\in B_D(x,r)$.
\end{proof}

Since the jump measure $\sB(x,y)|x-y|^{-d-\alpha}dxdy$ and the killing measure $\kappa(x)dx$ of  $(\sE^\kappa,  \sF^\kappa)$ are absolutely continuous with respect to  $m_d \otimes m_d$ and  $m_d$ respectively,   $X^\kappa$ satisfies the following L\'evy system formula: 
For any $x \in  D $, any non-negative Borel function $f$ on $D \times D_\partial$ vanishing on the diagonal, and any stopping time $\tau$,	
	\begin{align}\label{e:Levysystem-X-kappa}	
		\E_x\bigg[  \sum_{s \le \tau} f( X^\kappa_{s-},  X^\kappa_s)  \bigg] 
		= \E_x \left[ \int_0^\tau \bigg(  \int_{D} 
		\frac{f( X^\kappa_s, y) \sB( X^\kappa_s, y)}{| X^\kappa_s-y|^{d+\alpha}} dy  
		+ \kappa(X^\kappa_s)f(X_s^\kappa, \partial) \bigg) \,ds  \right].	
	\end{align}

For any open set $U\subset D$, we let $$\tau^\kappa_U:=\inf\{t>0: X^\kappa_t\notin U\}.$$ 
We  denote by $X^{\kappa,U}$ and $(P^{\kappa,U}_t)_{t \ge 0}$  the part of the process $X^\kappa$ killed upon exiting $U$ and its semigroup, respectively.
We denote   $(P^{\kappa,D}_t)_{t \ge 0}$ by $(P^{\kappa}_t)_{t \ge 0}$.  
Since $X^\kappa$ satisfies the interior parabolic H\"older regularity, for every open set $U\subset D$,  $(P^{\kappa,U}_t)_{t\ge 0}$ admits a jointly continuous heat kernel $p^\kappa_U(t,x,y)$ on $(0,\infty) \times U \times U$; see \cite[Theorem 4.13 and the paragraph below]{CKSV26a}.

The following  result was established  in \cite[Proposition 4.11]{CKSV26a}.
\begin{prop}\label{p:ndl2}
For any $R>0$ and $b\in (0,1)$, there exists $C=C(R,b)>0$ such that for any $x_0 \in D$, $r \in (0, (b\delta_D(x_0)) \wedge R]$ and $t\in (0, (br)^\alpha]$,
	it holds that
	\begin{equation*}
		p^{\kappa}_{B(x_0, r)}(t,z,y)  \ge  C t^{-d/\alpha}\quad  \text{for all } 
		y, z\in 	B(x_0, b t^{1/\alpha}).
	\end{equation*}
\end{prop}

The following lemma is a standard consequence of Proposition \ref{p:ndl2}; see, e.g., \cite[Proposition
3.5(ii)]{CKW20}. We omit the proof.

\begin{lemma}\label{l:E2-lower}
	For any $R>0$, there exist constants $C=C(R)>0$ and $\eps=\eps(R)\in (0,1)$ such that
	\begin{align}\label{e:E2-lower-1}
			\P_{x}(\tau^\kappa_{B(x,r)} \ge C r^\alpha) \ge \eps \quad \text{for all $x\in D$ and $r\in (0, \delta_D(x) \wedge R]$.} 
	\end{align}
	Consequently, 
	\begin{align*}
		\E_{x}[\tau^\kappa_{B(x,r)}] \ge 		C\eps r^\alpha \quad \text{for all $x\in D$ and $r\in (0, \delta_D(x) \wedge R]$.} 
	\end{align*}
\end{lemma}

\begin{lemma}\label{l:E2-upper}
 There exists  $C>0$ such that 	\begin{align*}
	\sup_{z\in B_D(x,r)}	\E_{z}[\tau^\kappa_{B_D(x,r)}] \le C r^\alpha \quad \text{for all $x\in D$ and $r\in (0, \wh R/8]$.} 
	\end{align*}
\end{lemma}
\begin{proof}
	Let $x\in D$ and $r\in (0,\wh R/8]$. We claim that there exist $w\in D$ and  a constant $c_1\in (0,1]$, depending only on $d, \wh R$ and $\Lambda$,  such that   $B(w,c_1r) \subset B_D(x,12r) \setminus B(x,2r)$. Indeed, if ,
	$\delta_D(x) \ge 4r$,
		then  any $w\in \partial B(x, 3r)$ satisfies $B(w,c_1r) \subset B_D(x,12r) \setminus B(x,2r)$ for any $c_1 \in (0,1]$.  
		Suppose $\delta_D(x)<4r$ and let $Q_x\in \partial D$ be such that $\delta_D(x)=|x-Q_x|$. Let $w:=7r\e_d$ in the coordinate system CS$_{Q_x}$. From the local coordinate system \eqref{e:local-coordinate_2}, we see that $\delta_D(w) \ge c_1r$ for some constant $c_1\in (0,1]$  depending only on $d,\wh R$ and $\Lambda$. Since $|w-x| \le |w-Q_x|+\delta_D(x)<11r$ and $|w-x| \ge |w-Q_x|-\delta_D(x)\ge 3r$, we have $B(w, c_1r) \subset B_D(x,12r) \setminus B(x,2r)$.   
		 By the L\'evy system formula \eqref{e:Levysystem-X-kappa} and \eqref{e:J-Phi}, since $\sB(x,y) \ge C_1^{-1}$, it follows that for all $z\in B_D(x,r)$,
	\begin{align*}
		1 &\ge \P_z \big( X^\kappa_{\tau^\kappa_{B(x,r)}} \in B(w,c_1r) \big)  = \E_z \left[ \int_0^{\tau^\kappa_{B(x,r)}} \int_{B(w,c_1r)} \frac{\sB(X^\kappa_s,y)}{|X^\kappa_s-y|^{d+\alpha}} dy ds\right]\\
		&\ge \frac{C_1^{-1} \E_z[\tau^\kappa_{B(x,r)}]}{(13r)^{d+\alpha}}\int_{B(w,c_1r)} dy  = c_2 r^{-\alpha}\E_z[\tau^\kappa_{B(x,r)}].
	\end{align*}
	This proves that $\E_z[\tau^\kappa_{B(x,r)}] \le c_2^{-1}r^\alpha$. 
\end{proof}

For any open set $U\subset D$, $x\in U$ and $R>0$, we define
$$
\Theta^U(x,R):=\E_x\left[\int_0^{\tau^\kappa_U}\Phi\left(\frac{R}{\delta_D(X^\kappa_s)}\right)\, ds \right].
$$

Recall that $p \in [(\alpha-1)_+, \alpha - \ub)$ is the unique constant satisfying \eqref{e:C(alpha,p,F)}.  
The next result is \cite[Theorem 6.2]{CKSV26a}. It is  one of the main results of \cite{CKSV26a},
and shows that the boundary decay of the respective quantities is of order $\delta_D(x)^p$. 
\begin{thm}\label{t:Dynkin-improve}
	There exist a constant $\eps_1\in (0, 1/8]$ and comparison constants such that  for any $Q\in \partial D$,    $R\in (0, \wh R]$ and $x\in B_D(Q,\eps_1R/2)$, 
	\begin{align*}
		&\P_x\big(X^\kappa_{\tau^\kappa_{B_D(Q,\eps_1R)}} \in D \big) \asymp R^{-\alpha}  \Theta^{B_D(Q,\eps_1R)}(x,R) \asymp  (\delta_D(x)/R)^p.
	\end{align*}
\end{thm}

 From Theorem \ref{t:Dynkin-improve}, we obtain the following estimate for the lifetime $\zeta^\kappa$; see \cite[Lemma 7.1]{CKSV26a} for the proof. 
\begin{lemma}\label{l:life}
	For any $T>0$,  there exists  $C=C(T) >0$ such that 
	$$
	\P_x (\zeta^\kappa>t) \le  C \left(1 \wedge \frac{\delta_D(x)}{t^{1/\alpha}}\right)^p \quad \text{for all $t\in (0,T]$ and $x\in D$}.
	$$
\end{lemma}

Using the  Markov property and 
Lemma \ref{l:life}, we obtain the following upper bounds for the  right-tail of $\tau^\kappa_{B_D(x,r)}$.

\begin{lemma}\label{l:survival-exponential}
There exist constants $C,C'>0$  such that for all $x\in D$, $r\in (0,\wh R/8]$ and 	$t\ge r^\alpha$, 
	\begin{align*}
		\P_x ( \tau^\kappa_{B_{D}(x,r)}>  t) \le  C \left( 1 \wedge \frac{\delta_D(x)}{r} \right)^pe^{-C' t/r^\alpha}.
	\end{align*}
\end{lemma}
\begin{proof} 
	Let $x\in D$ and $r\in (0,\wh R/8]$.  Set $B:=B_D(x,r)$.	By Lemma \ref{l:E2-upper},  there exists  $c_1>0$ independent of $x$ and $r$ such that 
	\begin{align}
		\label{e:nlbd1} 
		\sup_{z\in B} \E_z  \tau^\kappa_{B} \le 2^{-1} c_1  r^\alpha.
	\end{align}
	By Lemma \ref{l:life},  	since $t\mapsto	\P_x ( \tau^\kappa_{B}>t)$ is non-increasing, 
	it suffices to show that there exists a constant $c_2$ independent of $x$ and $r$ such that 
	\begin{align}\label{e:survival-exponential-claim}
		\P_x ( \tau^\kappa_{B}
		>  (c_1 n+1)  r^\alpha) \le c_2 \P_x ( \zeta^\kappa > r^\alpha) 
		2^{-n} \quad \text{for all $n \in \Z_+$}.
	\end{align}

	By Markov's inequality and \eqref{e:nlbd1}, 
	\begin{align}\label{e:nlbd2}
		&\sup_{z\in B}	\P_z ( 	\tau^\kappa_{B}
		> c_1r^\alpha)  \le  (c_1r^\alpha)^{-1} 	\sup_{z\in B} \E_z  \tau^\kappa_{B}  
		\le 2^{-1}.
	\end{align}
	Applying the Markov property at time $r^\alpha$ and at time $c_1r^\alpha$ ($n$ times) and then using  \eqref{e:nlbd2}, we get that for all $n \in \Z_+$, 
	\begin{align*}
		 \P_x ( \tau^\kappa_{B}
		> (c_1 n+1)r^\alpha)&\le  \E_x\left[	\P_{
			X^\kappa_{r^\alpha}} ( \tau^\kappa_{B}
		>c_1nr^\alpha) ;\,  \tau^\kappa_{B}
		>r^\alpha \right]\\
		&\le  \P_x ( \tau^\kappa_{B}
		>r^\alpha)  \sup_{z\in B}	\P_z ( 	\tau^\kappa_{B}	>c_1nr^\alpha) \\
		&\le \cdots \le  \P_x ( \tau^\kappa_{B}
		>r^\alpha) 	\Big( \sup_{z\in  B	}	\P_z ( 	\tau^\kappa_{B}
		>c_1r^\alpha)  \Big)^{n} \le 2^{-n} \P_x ( \zeta^\kappa>r^\alpha),
	\end{align*}where we used the fact that $\tau^\kappa_B\le \zeta^\kappa$ in the last inequality  
	Thus,  \eqref{e:survival-exponential-claim} holds and the proof is  complete. 
\end{proof}

The following proposition refines  Lemma \ref{l:E2-lower} for small $R$ by establishing lower bounds for the exit time up to the boundary.
\begin{prop}\label{p:met-bound-new} 
There exist constants $C>0$ and $\eps\in (0,1)$ 
 such that 
 	\begin{align}\label{e:met-bound-new-1} 
 	\P_{x}(\tau^\kappa_{B_D(x,r)
 	} \ge C r^\alpha) \ge \eps  \left( 1 \wedge \frac{\delta_D(x)}{r}\right)^p \quad \text{for all $r \in (0, \wh R]$ and $x\in D$.}
 \end{align}
 Consequently,
	\begin{align}\label{e:met-bound-new-2} 
		\E_x[\tau^\kappa_{B_D(x,r)} ]\ge C \eps \left( 1 \wedge \frac{\delta_D(x)}{r}\right)^p r^\alpha \quad \text{ for all $r \in (0, \wh R]$ and $x\in D$.}
	\end{align}
\end{prop}
\begin{proof}
It suffices to prove \eqref{e:met-bound-new-1}, since  \eqref{e:met-bound-new-2}  follows immediately from \eqref{e:met-bound-new-1} by Markov's inequality.

Let $r\in (0, \wh R]$, $x\in D$, and  $Q_x\in \partial D$ be such that $|x-Q_x|=\delta_D(x)$. Let $\eps_1\in (0,1/8]$ be the constant in Theorem \ref{t:Dynkin-improve}, and let  $c_1$ denote the constant $C$ in Lemma \ref{l:E2-lower} corresponding to $R=\wh R$. Let 
$c_2\in (0,c_1(\eps_1/2)^\alpha]$ be a constant to be chosen later.      If $\delta_D(x)\ge \eps_1 r/2$, then by 
Lemma \ref{l:E2-lower} (with $R=\wh R$), 
$$
	\P_{x}(\tau^\kappa_{B_D(x,r)} \ge c_2 r^\alpha) \ge 	\P_{x}(\tau^\kappa_{B(x,\eps_1 r/2)} \ge c_1 (\eps_1r/2)^\alpha) \ge c_3, 
$$
implying the desired inequality. 

Suppose $\delta_D(x)<\eps_1 r/2$. Since $D$ is $C^{1,\eta}$, there exist a constant $c_4\in (0,1/4)$, independent of $R$ and $x$, and $y\in B_D(Q_x, r/2)$ such that $B(y, 2c_4r) \subset B_D(Q_x, r/2) \setminus B(Q_x, \eps_1r)$. Note that $B(y,2c_4r) \subset B_D(x, r)$. 
Using the strong Markov property and the L\'evy system formula \eqref{e:Levysystem-X-kappa},  we obtain that
\begin{align}\label{e:mE-1}
& 	\P_{x}(\tau^\kappa_{B_D(x,r)} \ge c_2 r^\alpha) \ge 	\P_{x}\left(\tau^\kappa_{B_D(x,r)} \ge c_2 r^\alpha, \,  X^\kappa_{\tau^\kappa_{B_D(Q_x,\eps_1r)}}  \in B(y, c_4r)\right)\nn\\
	&\ge  \E_x\left[ \P_{X^\kappa_{\tau^\kappa_{B_D(Q_x,\eps_1r)}}}  (\tau^\kappa_{B_D(x,r)} \ge c_2 r^\alpha ); X^\kappa_{\tau^\kappa_{B_D(Q_x,\eps_1r)}}  \in B(y, c_4r)\right]\nn \\
	&\ge \left( \inf_{w\in B(y,c_4r)} \P_w (\tau^\kappa_{B_D(x,r)}  \ge c_2 r^\alpha ) \right) \E_x \left[ \int_0^{\tau^\kappa_{B_D(Q_x,\eps_1r)} } \int_{B(y,c_4r)} \frac{\sB(X^\kappa_s, w)}{|X^\kappa_s - w|^{d+\alpha}} dwds\right].
\end{align}
Note that $B(w, c_4r) \subset  B(y,2c_4r)  \subset B_D(x,r)$ for all $w \in B(y,c_4r)$. Hence, by  choosing $c_2:=(c_1(\eps_1/2)^\alpha) \wedge (c_1c_4^\alpha)$, 
 we get from \eqref{e:E2-lower-1} that  
\begin{align}\label{e:mE-2}
	\inf_{w\in B(y,c_4r)} \P_w (\tau^\kappa_{B_D(x,r)}  \ge c_2 r^\alpha ) \ge \inf_{w\in B(y,c_4r)} \P_w (\tau^\kappa_{B(w,c_4r)}  \ge 
	c_1(c_4r)^\alpha ) \ge c_3.
\end{align}
On the other hand, observe that $\delta_D(w) \asymp |z-w| \asymp r$ for all $z\in B_D(Q_x,\eps_1r)$ and $w \in B(y,c_4r)$.
Thus, by \eqref{e:J-Phi}, \eqref{e:Phi-scaling} (with $\beta=1$) and \eqref{e:without-A0}, we get that for all $z\in B_D(Q_x,\eps_1r)$ and $w \in B(y,c_4r)$,
\begin{align}\label{e:mE-3}
&	\frac{\sB(z,w)}{|z-w|^{d+\alpha}} \ge \frac{C_1^{-1}}{|z-w|^{d+\alpha}}\Phi \left(\frac{(|z-w|\wedge A_0)^2}{(\delta_D(z) \wedge A_0) (\delta_D(w) \wedge A_0)} \right) \ge \frac{c_5}{r^{d+\alpha}} \Phi \left(\frac{r\wedge A_0}{\delta_D(z) \wedge A_0} \right)\nn\\
	&\ge \frac{c_6}{r^{d+\alpha}} \left(1 \wedge \frac{A_0}{r}\right) \left( 1 \vee \frac{\delta_D(z)}{A_0} \right) \Phi \left(\frac{r}{\delta_D(z)} \right) \ge \frac{c_6}{r^{d+\alpha}} \left(1 \wedge \frac{A_0}{\wh R}\right) \Phi \left(\frac{r}{\delta_D(z)} \right).
\end{align}
Combining \eqref{e:mE-1} with \eqref{e:mE-2} and \eqref{e:mE-3}, and applying Theorem \ref{t:Dynkin-improve}, we arrive at
\begin{align*}
		\P_{x}(\tau^\kappa_{B_D(x,r)} \ge c_2 r^\alpha) &\ge  \frac{c_7m_d(B(y,c_4r))}{r^{d+\alpha}} \E_x \left[ \int_0^{\tau^\kappa_{B_D(Q_x,\eps_1r)} }  \Phi \left(\frac{r}{\delta_D(X^\kappa_s)} \right)   ds\right]  \\	&=c_8 r^{-\alpha} \Theta^{B_D(Q_x,\eps_1r)}(x,r) \ge c_9 (\delta_D(x)/r)^p .
\end{align*}

The proof is complete.
\end{proof}


\section{Small time heat kernel upper bound}\label{s:uHKE}

The goal of this section is to establish 
the sharp small time heat kernel upper bound. 
The main step is to prove the following bound, which involves only one boundary term. Combined with the semigroup property and Theorem \ref{t:CKSV25}, this estimate will yield the desired full upper bound. 
Recall that $\delta_D(x,t)$ is defined in \eqref{e:def-deltaDt}. 
\begin{prop}\label{p:blow-up-step2}
	There exists $C>0$ such that for all $t\in (0,  \wh{R}^{\alpha}]$ and $x,y \in D$,
	\begin{equation*}
		p^\kappa(t,x,y) \le C \left(1 \wedge 
		\frac{\delta_D(x) \wedge \delta_D(y)}{t^{1/\alpha}}\right)^p
		\left( t^{-d/\alpha} {\wedge} \frac{t}{|x-y|^{d+\alpha}}\Phi\left( \frac{(|x-y| \wedge A_0)^2}{(\delta_D(x,t) \wedge A_0)(\delta_D(y,t) \wedge A_0) }\right)\right).
	\end{equation*} 
\end{prop}

  To prove Proposition \ref{p:blow-up-step2}, we first establish several auxiliary results.  
 The proof of the proposition will be given later in this section. 
 We begin with the following lemma, which provides the desired estimate in the near-diagonal region. 
	\begin{lemma}\label{l:UHKD}
For any $T>0$,  there exists  $C=C(T) >0$ such that 
\begin{equation*}
	p^\kappa(t,x,y) \le C \left(1 \wedge \frac{\delta_D(x) \wedge \delta_D(y)}{t^{1/\alpha}}\right)^{p}t^{-d/\alpha} \quad \text{for all 
		$t\in (0,T]$ and $x,y \in D$}.
\end{equation*}
\end{lemma}
\begin{proof} By the semigroup property, \eqref{e:pk-le-pb} and Theorem \ref{t:CKSV25},  for all $t\in (0,T]$ and $x,y\in D$,
\begin{align*}		p^\kappa(t,x,y)& = \int_D p^\kappa(t/2,x,z)p^\kappa(t/2,z,y) dz
\le \int_D p^\kappa(t/2,x,z)\overline p(t/2,z,y) dz
\\
&\le c_1t^{-d/\alpha} \int_D p^\kappa(t/2,x,z)dz   = c_1t^{-d/\alpha} \P_x(\zeta^\kappa>t/2).\end{align*}
By symmetry, we also get $p^\kappa(t,x,y) \le  c_1t^{-d/\alpha} \P_y(\zeta^\kappa>t/2)$.
The result now follows from Lemma \ref{l:life}.
\end{proof}

To obtain  the desired off-diagonal upper heat kernel bound, we use the following lemma,  which follows from the strong Markov property and the L\'evy system formula by considering the space-time location at which the process exits a neighborhood of its initial point.  
For jump kernels without boundary blow-up, the unweighted version of this lemma is well known; 
see, e. g., \cite[Lemma 4.2]{CKS14} and  \cite[Lemma 2]{BGR10} and their proofs. 
However, in the present setting with boundary blow-up, when $ \overline{\beta} >\alpha/2$, both $\sup_{x,y \in D, \, |x-y|>\delta} p^\kappa(t,x,y)$ and $\sup_{x,y \in D, \, |x-y|>\delta} \overline p(t,x,y)$  may blow up as $t \to 0$ for every $\delta>0$. To account for this behavior, we introduce the following weighted version.

\begin{lemma}\label{l:general-upper-2}	
	Let  $V_1$ and $V_3$ be open subsets of $D$ with ${\rm dist}(V_1,V_3)>0$. Set $V_2:=D \setminus (V_1 \cup V_3)$. 
	For any 	positive real-valued Borel function $g$ on $(0, \infty)$,
	  $x\in V_1$, $y \in V_3$ and  $t>0$, it holds that 	
	\begin{align}
	&p^\kappa(t,x,y) \le    \P_x(\tau^\kappa_{V_1}<  t \wedge  \zeta^\kappa )
	\sup_{ 
		t/3 \le s\le t, \, z \in V_2} p^\kappa(s,z,y)
	\nn	\\	
		&\qquad \quad +  
		\P_x(\tau^\kappa_{V_1}>t/3) \bigg(\sup_{0<s\le t/3, \, z \in V_2  } 
		\frac{p^\kappa(s,z,y)}{g(s)} \bigg)	\left(\sup_{v \in V_1}\E_v \left[
			g( 2t/3-\tau^\kappa_{V_1});  	t/3 < \tau^\kappa_{V_1} <2 t/3 
	  \right]   \right) \nn\\	
		&\qquad \quad + {\rm dist}(V_1,V_3)^{-d-\alpha}\int_0^t\int_{V_3} \int_{V_1} p^{\kappa}_{V_1}(s, x, u) 
		\sB(u,w) p^\kappa(t-s, y,w) \,du dw ds.		\label{e:general-upper-2}
	\end{align} 
	\end{lemma} 
\begin{proof} 
	By the strong Markov property, the L\'evy system formula \eqref{e:Levysystem-X-kappa} and symmetry, we get that for any $x \in V_1$ and $y \in V_3$, 
	\begin{align*}	
	p^\kappa(t,x,y)	&=\E_x \left[p^\kappa(t-\tau^\kappa_{V_1}, X^\kappa_{\tau^\kappa_{V_1}}, y) ; \tau^{\kappa}_{V_1} < t \wedge \zeta^\kappa \right]  \\ 
	&	=\E_x \left[p^\kappa(t-\tau^\kappa_{V_1}, X^\kappa_{\tau^\kappa_{V_1}}, y) ; \tau^\kappa_{V_1}\le 2 t/3,
	\,  X^\kappa_{\tau^\kappa_{V_1}} \in V_2 \right]\\
		&\quad +
		\E_x \left[p^\kappa(t-\tau^\kappa_{V_1}, X^\kappa_{\tau^\kappa_{V_1}}, y) ; 2 t/3 <\tau^\kappa_{V_1}
		<		t,  X^\kappa_{\tau^\kappa_{V_1}} \in V_2\right]\\	
		&\quad +
		\E_x \left[p^\kappa(t-\tau^\kappa_{V_1}, X^\kappa_{\tau^\kappa_{V_1}}, y) ; \tau^\kappa_{V_1}<
		t,  X^\kappa_{\tau^\kappa_{V_1}} \in V_3 \right]\\	
		&=:I_1+I_2+I_3.		
	\end{align*}
We have 
$	I_1\le \P_x(\tau^\kappa_{V_1}<t \wedge  \zeta^\kappa)	\sup_{ t/3 \le  s \le t, \, z \in V_2} p^\kappa(s,z,y).$  For $I_3$, using the L\'evy system formula \eqref{e:Levysystem-X-kappa}, symmetry and the fact that $X^\kappa_s = X^{\kappa,V_1}_s$ for $s<\tau^\kappa_{V_1}$,  we get that
\begin{align*}
	I_3&= \E_x \left[ \int_0^{\tau^\kappa_{V_1} \wedge t} \int_{V_3} p^\kappa (t-s, y,w) \frac{\sB(X^\kappa_s, w)}{|X^\kappa_s-w|^{d+\alpha}} \, dwds  \right]\\
	&\le {\rm dist}(V_1,V_3)^{-d-\alpha} \E_x \left[ \int_0^{\tau^\kappa_{V_1} \wedge t} \int_{V_3} p^\kappa (t-s, y,w) \sB(X^{\kappa, V_1}_s, w) \, dwds   \right]\\
	&= {\rm dist}(V_1,V_3)^{-d-\alpha}\int_0^t\int_{V_3} \int_{V_1} p^{\kappa}_{V_1}(s, x, u) 
	\sB(u,w) p^\kappa(t-s, y,w) \,du dw ds.	
\end{align*}
For $I_2$, we note that  $\tau^\kappa_{V_1}= t/3+ \tau^\kappa_{V_1} \circ \theta^\kappa_{t/3}$ 
on $\{ \tau^\kappa_{V_1}>2t/3\}$
and $\zeta^\kappa= t/3+ \zeta^\kappa \circ \theta^\kappa_{t/3}$ 
on $\{\zeta^\kappa>2t/3\}$, where $\{\theta^\kappa_t:t\ge 0\}$ is the shift operator associated with $X^\kappa$. Hence,  by the Markov property, 
\begin{align*}	
I_2&=
\E_x \Big[\E_{X^{\kappa, V_1}_{t/3}}\big[p^\kappa(2t/3-\tau^\kappa_{V_1}, X^\kappa_{\tau^\kappa_{V_1}}, y) ; t/3 < \tau^\kappa_{V_1} 
<2 t/3, 
 X^\kappa_{\tau^\kappa_{V_1}} \in V_2\ \big]\Big] \\
&=
\E_x \bigg[\E_{X^{\kappa, V_1}_{t/3}}\bigg[ g(2 t/3-\tau^\kappa_{V_1}) \frac{p^\kappa(2t/3-\tau^\kappa_{V_1}, X^\kappa_{\tau^\kappa_{V_1}}, y)}{ g(2 t/3-\tau^\kappa_{V_1})} ;  t/3 < \tau^\kappa_{V_1} 
<2t/3,  X^\kappa_{\tau^\kappa_{V_1}} \in V_2\ \bigg]\bigg] \\ 
&\le \bigg( \sup_{0<s \le   t/3, \, z \in V_2  } \frac{p^\kappa(s,z,y)}{g(s)}  \bigg)
\int_{V_1} p^{\kappa}_{V_1} (t/3, x, v) \,\E_v \Big[g(2t/3-\tau^\kappa_{V_1}) ;  t/3 < \tau^\kappa_{V_1} <2 t/3 \Big]dv \\
&\le 	\P_x(\tau^\kappa_{V_1}>t/3) \bigg(\sup_{0<s \le   t/3, \, z \in V_2  } \frac{p^\kappa(s,z,y)}{g(s)} \bigg)	\left(\sup_{v \in V_1}\E_v \left[g(2  t/3-\tau^\kappa_{V_1});  
 t/3 < \tau^\kappa_{V_1} <2 t/3  \right]   \right).
\end{align*}
The proof is complete.
\end{proof}

The following  lemma will be used to estimate the first term  on the right-hand side of \eqref{e:general-upper-2}. 
\begin{lemma}\label{l:general-upper-first-term}
	For any $T>0$,	there exists $C=C(T)>0$ such that for all $t\in (0, T]$ and $x,y\in D$,
	\begin{align*}
		\sup_{t/3 \le s \le t, \, z \in D \setminus B(y, |x-y|/2)} \overline p(s,z,y)& \le  \frac{C t }{|x-y|^{d+\alpha}}
		\Phi\left( \frac{(|x-y| \wedge A_0)^2}{( t^{1/\alpha} \wedge A_0)(\delta_D(y,t) \wedge A_0)  }\right).
	\end{align*}
\end{lemma}
\begin{proof}
Let $t\in (0, T]$ and $x,y\in D$. 
	Using Theorem  \ref{t:CKSV25} in the first inequality below,  Lemma \ref{l:jump-density-monotonicty} in the second,
	and   \eqref{e:Phi-scaling} in the third, we obtain 
	\begin{align*}
		&	\sup_{t/3 \le s \le t, \, z \in D \setminus B(y, |x-y|/2)} \overline p(s,z,y) \nn\\
		&\qquad \le c_{1} 	\sup_{t/3 \le s \le t, \, z \in D \setminus B(y, |x-y|/2)} \frac{s}{|z-y|^{d+\alpha}} \Phi\left( \frac{(|z-y| \wedge A_0)^2}{(s^{1/\alpha} \wedge A_0)(\delta_D(y,s) \wedge A_0) }\right) \nn\\
		&\qquad \le \frac{c_{2} t }{|x-y|^{d+\alpha}}\sup_{ t/3 \le s \le t} 
		\Phi\left( \frac{(|x-y| \wedge A_0)^2}{( s^{1/\alpha} \wedge A_0)(\delta_D(y,s) \wedge A_0) }\right)\nn\\
		&\qquad =  \frac{c_{2} t }{|x-y|^{d+\alpha}}
		\Phi\left( \frac{(|x-y| \wedge A_0)^2}{( ( t/3)^{1/\alpha} \wedge A_0)(\delta_D(y, t/3) \wedge A_0)  }\right)\nn \\
		&\qquad\le \frac{c_{3} t }{|x-y|^{d+\alpha}}
		\Phi\left( \frac{(|x-y| \wedge A_0)^2}{( t^{1/\alpha} \wedge A_0)(\delta_D(y,t) \wedge A_0)  }\right).
	\end{align*}
	\end{proof}

	In the usual stable-like setting (namely, $\Phi=1$), the following lemma is straightforward, since the off-diagonal heat kernel estimates are of the form  $t |x-y|^{-d-\alpha}$, which is increasing in $t$. 	However, as mentioned above, in the present setting, the quantity  $\sup_{x,y \in D, \, |x-y|>\delta} \overline p(t,x,y)$  may blow up as $t \to 0$ for every $\delta>0$, and, in particular, need not be almost increasing in $t$.  By introducing a factor $t^\gamma$ for $\gamma>\ub/\alpha$, and replacing one boundary factor $\delta_D(x,t)$ by $t^{1/\alpha}$, we obtain the following lemma, which will be used to estimate the second term in \eqref{e:general-upper-2}.
	
	\begin{lemma}\label{l:general-upper-second-term}
		For any $\gamma >\ub/\alpha$ and $T>0$, there exists $C=C(\gamma,T)>0$ such that for  all $t \in (0, T]$ and $x,y \in D$,
		\begin{align*}
	 \sup_{0<s \le t, \, z \in D\setminus B(y,|x-y|/2)  } s^{\gamma} \overline p(s,z,y)  \le  \frac{C t^{1+\gamma} }{|x-y|^{d+\alpha}}
	 \Phi\left( \frac{(|x-y| \wedge A_0)^2}{(t^{1/\alpha} \wedge A_0)(\delta_D(y,t) \wedge A_0) }\right)	.
		\end{align*}
	\end{lemma}
	\begin{proof}
		Let $t\in (0,T]$ and $x,y \in D$.   Set $\beta:= (\ub + \alpha \gamma)/2$. Note that $\beta>\ub$ and $2\beta/\alpha<1+\gamma$.  
		 Using \eqref{e:Phi-scaling}  and \eqref{e:without-A0} twice, we obtain for all $s \in (0,t]$ and $z \in D\setminus B(y,|x-y|/2)$,
		\begin{align*}
			&s^{1+\gamma} \Phi\left( \frac{(|z-y| \wedge A_0)^2}{( \delta_D(z,s) \wedge A_0)(\delta_D(y,s) \wedge A_0) }\right)\le  s^{1+\gamma} \Phi\left( \frac{(|z-y| \wedge A_0)^2}{(s^{1/\alpha}  \wedge A_0)
				(\delta_D(y,s) \wedge A_0) }\right)\\
			&\le c_{1}s^{1+\gamma} \bigg(\frac{t^{1/\alpha}}{s^{1/\alpha}}\bigg)^{\beta}\bigg( \frac{\delta_D(y,t)}{\delta_D(y,s)}  \bigg)^{\beta}
			\Phi\left( \frac{(|z-y| \wedge A_0)^2}{(t^{1/\alpha} \wedge A_0)(\delta_D(y,t) \wedge A_0) }\right)\\
			&\le c_{1}  s^{1+\gamma-2\beta/\alpha} t^{2\beta/\alpha} 
			\Phi\left( \frac{(|z-y| \wedge A_0)^2}{(t^{1/\alpha} \wedge A_0)(\delta_D(y,t) \wedge A_0) }\right).
		\end{align*}
		Thus,  since $2\beta/\alpha<1+\gamma$,   Theorem \ref{t:CKSV25} and Lemma \ref{l:jump-density-monotonicty} imply that
		\begin{align*}
			&	 \sup_{0<s \le t, \, z \in D\setminus B(y,|x-y|/2)  } s^{\gamma} \overline p(s,z,y)  \nn\\
			&\qquad \le c_{2}  \sup_{0<s \le t, \, z \in D\setminus B(y,|x-y|/2)  } \frac{s^{1+\gamma}}{|z-y|^{d+\alpha}} 
			\Phi\left( \frac{(|z-y| \wedge A_0)^2}{(\delta_D(z,s) \wedge A_0)(\delta_D(y,s) \wedge A_0) }\right) \nn\\
			&\qquad\le c_{3} t^{2\beta/\alpha} 
		 \sup_{0<s\le t,\,  z \in D\setminus B(y,|x-y|/2)  } \frac{s^{1+\gamma - 2\beta/\alpha}}{|z-y|^{d+\alpha}}
			\Phi\left( \frac{(|z-y| \wedge A_0)^2}{(t^{1/\alpha} \wedge A_0)(\delta_D(y,t) \wedge A_0) }\right)\nn	\\ 
			&\qquad\le  \frac{c_{4} t^{1+\gamma} }{|x-y|^{d+\alpha}}
			\Phi\left( \frac{(|x-y| \wedge A_0)^2}{(t^{1/\alpha} \wedge A_0)(\delta_D(y,t) \wedge A_0) }\right)	.
		\end{align*}
	\end{proof}

For the second term on the right-hand side of \eqref{e:general-upper-2}, we will choose $g$  to be a negative power function. 
To estimate this term, we use the following estimates for the time derivatives of Dirichlet heat kernels, which follow from the ultracontractivity of    $(P^{\kappa, U}_t)_{t\ge 0}$. The argument is well known (see \cite[Corollary 5.7]{GT} or \cite[(4.1.4)]{S}), but we include it here for completeness.

\begin{prop}\label{p:time-derivative}
There exists  $C>0$ such that for any open set $U\subset D$ and $t>0$,
	\begin{align*}
		\sup_{x,y \in U} |\partial_t 	p^{\kappa}_U(t,x,y)| \le Ct^{-1-d/\alpha}e^{t/ \wh R^\alpha} \, .	
	\end{align*}
\end{prop}
\begin{proof}
Let $U \subset D$ be an open subset. By the ultracontractivity of the semigroups of $\overline X$ in \cite[Proposition 3.12]{CKSV25} and the continuity of heat kernels, we have
	\begin{align}\label{e:Ultracontractivity}\sup_{x,y \in U} |	p^{\kappa}_U(t,x,y)| \le \sup_{x,y \in U} |	\overline p(t,x,y)|  \le  c_1 t^{-d/\alpha}e^{t/ \wh R^\alpha} \quad \text{for all $t>0$}.	\end{align}
	Moreover, by \cite[Corollary 5.5 and Lemma 5.6]{GT}, since $x\mapsto p^\kappa_U(t,x,y)$ is continuous on $U$, the function $t\mapsto p^\kappa_U(t,x,y)$ is differentiable on $(0,\infty)$ for all  $x,y \in U$, and satisfies that
	\begin{gather}
		\partial_t p^\kappa_U(t,x,y)  = \int_U p^\kappa_U(s,x,z) \partial_t p^\kappa_U(t-s, z,y) dz,\label{e:time-derivative-1}\\
			\int_U \left( \partial_t p^\kappa_U(t,z,y) \right) ^2 dz \le \frac{ p^\kappa_U(2s,y,y)}{(t-s)^2} ,\label{e:time-derivative-2}
	\end{gather} 
	for all $0<s<t$ and $x,y \in U$. Combining \eqref{e:time-derivative-1}, \eqref{e:time-derivative-2} and \eqref{e:Ultracontractivity}, and using the symmetry and the semigroup property,  we obtain that for all $t>0$ and $x,y \in U$,
	\begin{align*}
		& |\partial_t 	p^{\kappa}_U(t,x,y)|   = \left|\int_U p^\kappa_U(t/4,x,z) \partial_t p^\kappa_U(3t/4, z,y) dz \right| \\
		&\le \left( \int_U  p^\kappa_U(t/4,x,z)p^\kappa_U(t/4,z,x) dz\right)^{1/2} \left( \int_U  \partial_t p^\kappa_U(3t/4,z,y)^2 dz\right)^{1/2} \\
		&\le \frac{\left( p^\kappa_U(t/2,x,x) p^\kappa_U(t/2,y,y) \right)^{1/2}}{t/2} \le c_2 t^{-1-d/\alpha} e^{t/\wh R^\alpha}.
	\end{align*}
\end{proof}

\begin{prop}
	\label{p:exitu}
	For every $\gamma, \eps  \in (0,1)$ and $T>0$, there exists $C=C(\gamma,\eps, T)$ such that for any open subset $V$ of $D$,
	$$
\sup_{v\in V}	\E_v \left[(t-\tau^\kappa_{V})^{-\gamma};  \eps t <\tau^\kappa_{V}< t  \right]  
	\le C m_d(V)t^{-d/\alpha-\gamma} \quad \text{for all $t\in (0,T]$}.
	$$
\end{prop}
\begin{proof}
Let $V\subset D$ be an open subset and $t\in (0,T]$.	Using the identity 
	$$
	\gamma	\int_{t-a}^{t} s^{-1-\gamma} ds=(t-a)^{-\gamma}-t^{-\gamma},  \quad a \in (0, t), $$
	we get that for all $v \in V$, 
	\begin{align}		\label{e:exitu-1}
		& \E_v \left[(t-\tau^\kappa_{V})^{-\gamma};  \eps t <\tau^\kappa_{V}< t  \right]   \nn\\
		&=t^{-\gamma} \P_v \left( \eps t <\tau^\kappa_{V}< t  \right)+\gamma \, \E_v \bigg[
		\int_{t-\tau^\kappa_{V}}^{t} s^{-1-\gamma} ds
		; \, \eps t <\tau^\kappa_{V}< t  \bigg]   \nn\\
		&= t^{-\gamma} \P_v \left( \eps t <\tau^\kappa_{V}< t  \right)+\gamma  \int_0^{t} \P_v \big( (t-s) \vee \eps t < \tau^\kappa_{V}< t  \big) s^{-1-\gamma}  ds\, . 
	\end{align}
	Note that
	\begin{align*}	
		&\P_v \left( (t-s) \vee (\eps t) <\tau^\kappa_{V} < t  \right)\\
		&\le 
		\P_v \left( \tau^\kappa_{V}>(t-s) \vee (\eps t)  \right) - \P_v \left(\tau^\kappa_{V} > t  \right)= \int_V \left(p^{\kappa}_V((t-s) \vee (\eps t), v, w) - p^{\kappa}_V(t, v, w)\right)  dw \\
		&\le m_d(V)\sup_{z, w \in V} \left|p^{\kappa}_V((t-s) \vee (\eps t), z, w) - p^{\kappa}_V(t, z, w)\right|.
	\end{align*}
	By the mean value inequality and  Proposition \ref{p:time-derivative}, for all $0<s\le t$,		
	\begin{align*}	
		&\sup_{z, w \in V} \left|p^{\kappa}_V((t-s) \vee (\eps t), z, w) - p^{\kappa}_V(t, z, w)\right|\le
		s\Big(\sup_{a \in [\eps t, t], z, w \in V} 	|\partial_t p^{\kappa}_V(a, z, w)| \Big)  
		\le	 c_1 st^{-d/\alpha-1}. 
	\end{align*}
	Thus, we obtain that for all $0<s\le t$,
	\begin{align}	
		\label{e:exitu-2}
		\P_v \left( (t-s) \vee (\eps t) <\tau^\kappa_{V} < t  \right)\le c_1 m_d(V)st^{-d/\alpha-1}. 
	\end{align}
 Applying \eqref	{e:exitu-2} twice (with $s$ and with $s$ replaced by $t$), and using the fact that $\gamma<1$,  we conclude from \eqref{e:exitu-1} that  
	\begin{align*}	
		\E_v \left[(t-\tau^\kappa_{V})^{-\gamma};  \eps t <\tau^\kappa_{V} < t \right]  
	&	\le c_1 m_d(V) \left(  t^{-d/\alpha-\gamma} + \gamma t^{-d/\alpha-1}\int_0^t s^{-\gamma}ds  \right) = c_2 m_d(V) t^{-d/\alpha-\gamma}.
	\end{align*}
\end{proof}

To estimate the third term on the right-hand side of  \eqref{e:general-upper-2}, we will use suitable weighted integral estimates for the heat kernels to account for the blow-up behavior of $\sB(u,w)$ near the boundary.
We first recall the following result for the (untruncated) process $\overline X$ from \cite[Proposition 5.12]{CKSV25}. 
 
	 \begin{prop}\label{p:blow-up-integrated-improve-final}	There exists $C>0$ such that for all $a>0$, $x\in D$  
	 and $t\in (0, \wh R^\alpha]$,
	 	\begin{align*}
	 		\int_D \overline p (t, x, z) \Phi \bigg( \frac{a}{\delta_D(z)\wedge A_0}\bigg)   dz\,\le C\Phi \bigg( \frac{a}{\delta_D(x,t)\wedge A_0}\bigg) .
	 	\end{align*}
	 \end{prop}

	 Using the above proposition, we obtain
	 	 \begin{lemma}\label{l:blow-up-integrate-1}
	 	 	There exists $C>0$ such that for all $a>0$, $x\in D$ and  $t\in (0, \wh R^\alpha]$,
	 	\begin{align*}
	 	\int_0^t	\int_D p^\kappa (s, x, u) \Phi \bigg( \frac{a}{\delta_D(u)\wedge A_0}\bigg)   duds \,
	 	\le C t	\left( 1 \wedge \frac{\delta_D(x)}{t^{1/\alpha}}\right)^p \Phi \bigg( \frac{a}{t^{1/\alpha}\wedge A_0}\bigg).
	 	\end{align*}
	 \end{lemma}
	 \begin{proof}
	 Let $a>0$, $x\in D$ and $t\in (0, \wh R^\alpha]$.	Choose $\beta \in (\overline{\beta}, 1\wedge \alpha)$ such that  $\beta <\alpha-p$.
	Note that, by Proposition \ref{p:blow-up-integrated-improve-final}, \eqref{e:Phi-scaling}
	 and \eqref{e:without-A0}, 
	 \begin{align*}
	 	&\int_0^t	\int_D p^\kappa (s, x, u) \Phi \bigg( \frac{a}{\delta_D(u)\wedge A_0}\bigg)   duds 
		\le \int_0^t	\int_D \overline p (s, x, u) \Phi \bigg( \frac{a}{\delta_D(u)\wedge A_0}\bigg)   duds\nn\\
		 &\le c_1 \int_0^t	 \Phi \bigg( \frac{a}{s^{1/\alpha}\wedge A_0}\bigg)ds
		\le c_2 t^{\beta/\alpha} \Phi \bigg( \frac{a}{t^{1/\alpha}\wedge A_0}\bigg) \int_0^t	s^{-\beta/\alpha}  ds 
	    =c_3 t \Phi \bigg( \frac{a}{t^{1/\alpha}\wedge A_0}\bigg).
	 	\end{align*}
		Thus, when $\delta_D(x) \ge t^{1/\alpha}$, the desired inequality holds.
		
		Suppose $\delta_D(x) < t^{1/\alpha}$. By the semigroup property,  Proposition \ref{p:blow-up-integrated-improve-final}, \eqref{e:Phi-scaling} and \eqref{e:without-A0}, 
	 \begin{align}\label{e:blow-upn2}
	 &\int_0^t	\int_D p^\kappa (s, x, u) \Phi \bigg( \frac{a}{\delta_D(u)\wedge A_0}\bigg)   duds \nn\\
	 	&=\int_0^t	\int_D \int_D p^\kappa (s/2, x, z) p^\kappa (s/2, z, u) \Phi \bigg( \frac{a}{\delta_D(u)\wedge A_0}\bigg)   dudzds \nn\\
		&\le \int_0^t	\int_Dp^\kappa (s/2, x, z) \left( \int_D  \overline p (s/2, z, u) \Phi \bigg( \frac{a}{\delta_D(u)\wedge A_0}\bigg)   du
		\right)dzds \nn\\
		&\le c_4  \int_0^t 	
	\int_D p^\kappa (s/2, x, z) dz \,
		 \Phi \bigg( \frac{a}{s^{1/\alpha}\wedge A_0}\bigg) ds \nn\\
		&\le c_5 t^{\beta/\alpha} \Phi \bigg( \frac{a}{t^{1/\alpha}\wedge A_0}\bigg) \int_0^t 	s^{-\beta/\alpha} 
	\int_D p^\kappa (s/2, x, z) dz\, ds.
	 	\end{align}
		Thus,  using Lemma \ref{l:life} and the fact that $p+\beta<\alpha$, we obtain
			from \eqref{e:blow-upn2} that 
			\begin{align*}
					&\int_0^t	\int_D p^\kappa (s, x, u) \Phi \bigg( \frac{a}{\delta_D(u)\wedge A_0}\bigg)   duds\\
					& \le c_6 t^{\beta/\alpha} \delta_D(x)^p\Phi \bigg( \frac{a}{t^{1/\alpha}\wedge A_0}\bigg) \int_0^t 	s^{-(\beta+p)/\alpha}ds = 
					c_7 t	\left(  \frac{\delta_D(x)}{t^{1/\alpha}}\right)^p \Phi \bigg( \frac{a}{t^{1/\alpha}\wedge A_0}\bigg)
				.
			\end{align*}
			The proof is complete.
		\end{proof}	  

\begin{lemma}\label{l:blow-up-integrate-2} 
  For any $\gamma \in [0,p+1)$ and $T>0$, there exists  $C=C(\gamma,T) >0$ such that, for  all $t\in (0,T]$, $s\in [t/2,t]$,  $x_0 \in \overline D$ and  $x\in D$,
   \begin{align*} 
\int_{B_D(x_0, t^{1/\alpha})} p^{\kappa}(s, x, u) \delta_D(u)^{-\gamma}  du  
  \le  C	t^{-\gamma/\alpha} \left( 1 \wedge \frac{\delta_D(x)}{t^{1/\alpha}}\right)^p.
	\end{align*}
\end{lemma}
\begin{proof}
Let  $t\in (0,T]$, $s\in [t/2,t]$ and $x_0 \in \overline D$. Set $V:=B_D(x_0,t^{1/\alpha})$. 
By Lemma \ref{l:UHKD}, there exists $c_1=c_1(T)>0$   such that 
$$
\sup_{v,w \in V} p^{\kappa}(s/3, v, w) \le \sup_{v,w \in D} p^{\kappa}(s/3, v, w)  \le c_1 (t/6)^{-d/\alpha}.
$$
Using this,  the symmetry  and semigroup property of $p^\kappa$,  and Lemma \ref{l:life},  we obtain
for all $x \in D$, 
\begin{align*}
	 	&\int_{V} p^{\kappa}(s, x, u) \delta_D(u)^{-\gamma}  du  \\
		&=\int_{V}\int_{D}\int_{D}
			p^{\kappa}(s/3, x, v)p^{\kappa}(s/3, v, w)p^{\kappa}(s/3, u,w )dvdw\, \delta_D(u)^{-\gamma} du  \\
		& \le c_2 t^{-d/\alpha}\left( \int_{D}p^{\kappa}(s/3, x, v) dv \right)\int_{V}\left(\int_{D}
			p^{\kappa}(s/3, u, w)dw \right) \delta_D(u)^{-\gamma} du  \\
			&\le  c_3 t^{-d/\alpha}\left(1 \wedge  \frac{\delta_D(x)}{(s/3)^{1/\alpha}}\right)^p \int_{V} \left(1 \wedge  \frac{\delta_D(u)}{(s/3)^{1/\alpha}}\right)^p \delta_D(u)^{-\gamma} du\\
		& \le   6^{2p/\alpha}c_3 t^{-d/\alpha} 	 \left(1 \wedge  \frac{\delta_D(x)}{t^{1/\alpha}}\right)^p \int_{V}   \left(1 \wedge  \frac{\delta_D(u)}{t^{1/\alpha}}\right)^p \delta_D(u)^{-\gamma} du.
\end{align*}
Thus, it suffices to show that
\begin{align}\label{e:blow-up-integrate-2-claim}
 I:=\int_{V}   \left(1 \wedge  \frac{\delta_D(u)}{t^{1/\alpha}}\right)^p \delta_D(u)^{-\gamma} du\le c_4 t^{(d-\gamma)/\alpha}
\end{align}
for some constant $c_4>0$ independent of $t$ and $x_0$. If $\delta_D(x_0) \ge 2t^{1/\alpha}$, then $\delta_D(u) >t^{1/\alpha}$ for all $u \in V$. Hence, 
$$
I \le t^{-\gamma/\alpha} m_d(B(x_0,t^{1/\alpha})) = c_5 t^{(d-\gamma)/\alpha}.
$$
 If $\delta_D(x_0)<2t^{1/\alpha}$, then since $p-\gamma>-1$, we get from  Lemma \ref{l:Aikawa-integrals} that 
$$
 I\le  t^{-p/\alpha}\int_{V}  \delta_D(u)^{p-\gamma} du  \le c_6 t^{(d-p)/\alpha} ( \delta_D(x_0) \vee t^{1/\alpha})^{p-\gamma} \le c_7 t^{(d-\gamma)/\alpha}.$$
This completes the proof of \eqref{e:blow-up-integrate-2-claim}, and thus, of the lemma. 
\end{proof}

 We now give the proof of Proposition \ref{p:blow-up-step2}.
	 
\begin{proof}[Proof of Proposition \ref{p:blow-up-step2}]
	 Let $t\in (0, \wh R^\alpha]$ and $x,y \in D$. By symmetry, we may assume without loss of generality that  $\delta_D(x) \le \delta_D(y)$.
If 
$\delta_D(x) \ge 2^{-4}t^{1/\alpha}$ or $|x-y| \le 8t^{1/\alpha}$, then the result follows from  
Theorem \ref{t:CKSV25} combined with \eqref{e:pk-le-pb} or Lemma \ref{l:UHKD}, respectively.

Assume now that  
$\delta_D(x)<2^{-4}t^{1/\alpha}$ and $|x-y|>8t^{1/\alpha}$. For the desired result, it suffices to show that there exists $c_1>0$ independent of $t,x,y$ such that
\begin{equation}\label{e:UHK-claim}
p^\kappa(t,x,y)\le 	\frac{c_1 \delta_D(x)^p t^{1-p/\alpha}}{|x-y|^{d+\alpha}}\Phi\left( \frac{(|x-y| \wedge A_0)^2}{(t^{1/\alpha} \wedge A_0)(\delta_D(y,t) \wedge A_0)}\right).
\end{equation}
To prove \eqref{e:UHK-claim}, we will use Lemma \ref{l:general-upper-2}.

Let $\beta \in (\ub, 1\wedge \alpha)$ and let $Q_x\in \partial D$ satisfy $|x-Q_x|=\delta_D(x)$. Set $V_1:=B_D(Q_x, t^{1/\alpha})$,  
$V_3:=B_D(y, |x-y|/2)$ and $V_2:=D \setminus (V_1 \cup V_3)$. Applying Lemma \ref{l:general-upper-2} with $g(s)=s^{-\beta/\alpha}$, we obtain that 
\begin{align*}
	p^\kappa(t,x,y)\le I_1+I_2+I_3,
\end{align*}
where 
\begin{align*}
	I_1&:=
	\P_x(\tau^\kappa_{V_1}<t\wedge  \zeta^\kappa)
	\sup_{ t/3 \le s \le t, \, z \in V_2} p^\kappa(s,z,y),\\
	I_2&:=	\P_x(\tau^\kappa_{V_1}>t/3) \bigg(\sup_{0<s \le   t/3, \, z \in V_2  } s^{\beta/\alpha}p^\kappa(s,z,y) \bigg)	\sup_{v \in V_1}\E_v \left[( 2t/3-\tau^\kappa_{V_1})^{-\beta/\alpha};  
	 \tau^\kappa_{V_1} \in (t/3,2t/3)
	\right]   ,
\end{align*}
and
\begin{align*}
	I_3:= {\rm dist}(V_1,V_3)^{-d-\alpha}\int_0^t\int_{V_3} \int_{V_1} p^{\kappa}_{V_1}(s, x, u) 
	\sB(u,w) p^\kappa(t-s, y,w) \,du dw ds.	
\end{align*}

For $I_1$, we first show that
\begin{align}\label{e:UHK-claim-I1-1}
	\P_x(\tau^\kappa_{V_1}<t\wedge  \zeta^\kappa) \le c_2(\delta_D(x)/t^{1/\alpha})^p,
\end{align}
for some constant $c_2$ independent of $t$ and $x$. 
Indeed, if $\delta_D(x)\ge \eps_1t^{1/\alpha}/2$, where  $\eps_1\in (0,1/8]$ is the constant in Theorem \ref{t:Dynkin-improve}, then  \eqref{e:UHK-claim-I1-1} holds trivially by taking $c_2$ larger than $(\eps_1/2)^{-p}$. If $\delta_D(x)< \eps_1t^{1/\alpha}/2$, then by Theorem \ref{t:Dynkin-improve}, we get
\begin{align*}
	\P_x(\tau^\kappa_{V_1}<t\wedge  \zeta^\kappa) \le \P_x\big(X^\kappa_{\tau^\kappa_{B_D(Q_x, \eps_1t^{1/\alpha})}} \in D\big) \le c_3(\delta_D(x)/t^{1/\alpha})^p,
\end{align*} 
completing the proof of \eqref{e:UHK-claim-I1-1}. Combining \eqref{e:UHK-claim-I1-1} with Lemma \ref{l:general-upper-first-term},  
we obtain
\begin{align*}
	I_1& \le c_2 \delta_D(x)^p t^{-p/\alpha}
	\sup_{ t/3 \le s \le t, \, z \in D\setminus B(y, |x-y|/2)} \overline p(s,z,y) \\
	&\le  \frac{c_4\delta_D(x)^p t^{1-p/\alpha} }{|x-y|^{d+\alpha}}
	\Phi\left( \frac{(|x-y| \wedge A_0)^2}{( t^{1/\alpha} \wedge A_0)(\delta_D(y,t) \wedge A_0)  }\right).
\end{align*}
For $I_2$, combining Lemma \ref{l:life}, Lemma \ref{l:general-upper-second-term} and Proposition \ref{p:exitu}, we obtain that
\begin{align*}
	I_2&\le \P_x(\zeta^\kappa>t/3) \bigg(\sup_{0<s \le   t, \, z \in D \setminus B(y,|x-y|/2)  } s^{\beta/\alpha}\overline p(s,z,y) \bigg)	\\
	&\quad \times \sup_{v \in V_1}\E_v \left[( 2t/3-\tau^\kappa_{V_1})^{-\beta/\alpha};  
	\tau^\kappa_{V_1} \in (t/3,2t/3)
	\right]  \\
	&\le c_5 \left( \frac{\delta_D(x)}{t^{1/\alpha}}\right)^p\frac{ m_d(V_1) t^{1+\beta/\alpha -d/\alpha-\beta/\alpha} }{|x-y|^{d+\alpha}}
	\Phi\left( \frac{(|x-y| \wedge A_0)^2}{(t^{1/\alpha} \wedge A_0)(\delta_D(y,t) \wedge A_0) }\right)\\
	&\le \frac{ c_3\delta_D(x)^p t^{1-p/\alpha} }{|x-y|^{d+\alpha}}
	\Phi\left( \frac{(|x-y| \wedge A_0)^2}{(t^{1/\alpha} \wedge A_0)(\delta_D(y,t) \wedge A_0) }\right).
\end{align*}

Now, we estimate $I_3$. By the 
triangle inequality, for all $u \in V_1$ and $w \in V_3$, we have
\begin{equation*}
|u-w| \ge |x-y| - |x-u| - |y-w| \ge |x-y| - 2t^{1/\alpha} - \frac{|x-y|}{2} \ge \frac{|x-y|}{4} >2t^{1/\alpha},
\end{equation*}
and
\begin{equation*}
 |u-w| \le |x-y| + |x-u| + |y-w| \le |x-y| + 2t^{1/\alpha} + \frac{|x-y|}{2} \le 2 |x-y| . 
\end{equation*}
Hence, dist$(V_1,V_3) \ge |x-y|/4$ and, by \eqref{e:J-Phi}, \eqref{e:Phi-scaling} and \eqref{e:without-A0},
\begin{equation*}
	\sB(u,w) \le
	c_6 \Phi\left(\frac{(|x-y| \wedge A_0)^2}{(\delta_D(u)\wedge A_0) (\delta_D(w)\wedge A_0)}\right)\, \quad \text{for all $u \in V_1$ and $w\in V_3$}.
\end{equation*}
Therefore,
\begin{align*}
	I_3&\le \frac{c_7}{|x-y|^{d+\alpha}}\int_0^{t/2}  \int_{V_1} p^{\kappa}(s, x, u) \,
\bigg[ \int_{V_3} p^\kappa(t-s, y,w) \Phi\left(\frac{(|x-y| \wedge A_0)^2}{(\delta_D(u)\wedge A_0) (\delta_D(w)\wedge A_0)}\right)  dw \bigg] \, du ds\\
&\quad + \frac{c_7}{|x-y|^{d+\alpha}}\int_{t/2}^t \int_{V_1} p^{\kappa}(s, x, u)  \bigg[ \int_{V_3} 
 p^\kappa(t-s, y,w) \Phi\left(\frac{(|x-y| \wedge A_0)^2}{(\delta_D(u)\wedge A_0) (\delta_D(w)\wedge A_0)}\right) dw \bigg]   du  ds	\\
&=:I_{3,1}+I_{3,2}.
\end{align*}

Applying Proposition \ref{p:blow-up-integrated-improve-final}  (with $a=(|x-y|\wedge A_0)^2/(\delta_D(u)\wedge A_0)$),  we obtain that for all $s\in (0,t)$ and  $u \in V_1$,
\begin{align}\label{e:UHK-1}
 &\int_{V_3}   p^\kappa(t-s, y,w) \Phi\left(\frac{(|x-y| \wedge A_0)^2}{(\delta_D(u)\wedge A_0) (\delta_D(w)\wedge A_0)}\right) dw 
 \nn\\
 &\le   \int_{D}   \overline p(t-s, y,w) 
 \Phi\left(\frac{(|x-y| \wedge A_0)^2}{(\delta_D(u)\wedge A_0) (\delta_D(w)\wedge A_0)}\right) dw\nn\\ 
 &\le c_8
 \Phi\left(\frac{(|x-y| \wedge A_0)^2}{(\delta_D(u)\wedge A_0) (\delta_D(y,t-s)\wedge A_0)}\right).
\end{align}
Note that $\delta_D(u) \le |u-Q_x| <t^{1/\alpha}$ for all $u \in V_1$. Thus, using  \eqref{e:Phi-scaling} and \eqref{e:without-A0}, we deduce from \eqref{e:UHK-1} that for all $s\in (0,t/2]$ and $u \in V_1$,
\begin{align}\label{e:UHK-1-small}
	&\int_{V_3}   p^\kappa(t-s, y,w) \Phi\left(\frac{(|x-y| \wedge A_0)^2}{(\delta_D(u)\wedge A_0) (\delta_D(w)\wedge A_0)}\right) dw\nn\\
	& \le  c_9
	\Phi\left(\frac{(|x-y| \wedge A_0)^2}{(\delta_D(u)\wedge A_0) (\delta_D(y,t)\wedge A_0)}\right),
\end{align}
and that for all $s\in (t/2,t)$ and $u \in V_1$,
\begin{align}\label{e:UHK-1-large}
	&	\int_{V_3}   p^\kappa(t-s, y,w) \Phi\left(\frac{(|x-y| \wedge A_0)^2}{(\delta_D(u)\wedge A_0) (\delta_D(w)\wedge A_0)}\right) dw\nn\\
	&\le c_{10}\bigg(\frac {t^{1/\alpha}} {\delta_D(u)}\bigg)^{\beta} \bigg(\frac {\delta_D(y,t)} {\delta_D(y,t-s)}\bigg)^{\beta} \Phi\left(\frac{(|x-y| \wedge A_0)^2}{(t^{1/\alpha}\wedge A_0)  (\delta_D(y,t)\wedge A_0)}\right) \nn\\ 
	&\le 
	\frac{c_{10}t^{2\beta/\alpha}}{\delta_D(u)^\beta (t-s)^{\beta/\alpha}} \Phi\left(\frac{(|x-y| \wedge A_0)^2}{(t^{1/\alpha}\wedge A_0)  (\delta_D(y,t)\wedge A_0)}\right).
\end{align}
Using \eqref{e:UHK-1-small} and applying  Lemma \ref{l:blow-up-integrate-1} with $a=(|x-y| \wedge A_0)^2/(\delta_D(y,t)\wedge A_0)$, we  get 
\begin{align}
I_{3,1}	&\le \frac{c_{11}}{|x-y|^{d+\alpha}} \int_0^{ t }\int_{V_1} p^{\kappa}(s, x, u) \, 	\Phi\left(\frac{(|x-y| \wedge A_0)^2}{(\delta_D(u)\wedge A_0) (\delta_D(y,t)\wedge A_0)}\right) du ds \nn\\
	& \le\frac{ c_{12}\delta_D(x)^p t^{1-p/\alpha} }{|x-y|^{d+\alpha}}
	\Phi\left( \frac{(|x-y| \wedge A_0)^2}{(t^{1/\alpha} \wedge A_0)(\delta_D(y,t) \wedge A_0) }\right).
\end{align}
On the other hand,  using \eqref{e:UHK-1-large} and applying Lemma \ref{l:blow-up-integrate-2}, since 
$\beta< \alpha \wedge (p+1)$, we obtain
\begin{align}
I_{3,2}	& \le
 \frac{c_{13}t^{2\beta/\alpha}}{|x-y|^{d+\alpha}}\Phi\left(\frac{(|x-y| \wedge A_0)^2}{(t^{1/\alpha}\wedge A_0)  (\delta_D(y,t)\wedge A_0)}\right) \int_{t/2}^t (t-s)^{-\beta/\alpha}  \int_{V_1} p^\kappa (s,x,u) \delta_D(u)^{-\beta}   du \, ds\\
	& \le \frac{c_{14}t^{2\beta/\alpha}}{|x-y|^{d+\alpha}} \Phi\left(\frac{(|x-y| \wedge A_0)^2}{(t^{1/\alpha}\wedge A_0)  (\delta_D(y,t)\wedge A_0)}\right)  \int_{t/2}^t (t-s)^{-\beta/\alpha}  s^{-\beta/\alpha} \bigg( \frac{\delta_D(x)}{s^{1/\alpha}}\bigg)^{p} \, ds\\
	& \le \frac{c_{14}\delta_D(x)^pt^{2\beta/\alpha}}{(t/2)^{(\beta+p)/\alpha}|x-y|^{d+\alpha}} \Phi\left(\frac{(|x-y| \wedge A_0)^2}{(t^{1/\alpha}\wedge A_0)  (\delta_D(y,t)\wedge A_0)}\right)  \int_{t/2}^t (t-s)^{-\beta/\alpha}   ds\\
	&= \frac{ c_{15}\delta_D(x)^p t^{1-p/\alpha} }{|x-y|^{d+\alpha}}
	\Phi\left( \frac{(|x-y| \wedge A_0)^2}{(t^{1/\alpha} \wedge A_0)(\delta_D(y,t) \wedge A_0) }\right).
\end{align}

The proof is complete.
\end{proof}

We are ready to establish the desired full upper bound. 

Define for $t>0$ and $x,y \in D$,
\begin{align}\label{e:def-wtp}
	&	\wt p(t, x,y):= 
	t^{-d/\alpha} \wedge  \left( \frac{ t }{|x-y|^{d+\alpha}}\Phi\left( \frac{(|x-y| \wedge A_0)^2}{(\delta_D(x,t) \wedge A_0)(\delta_D(y,t) \wedge A_0) }\right) \right). 
\end{align}

\begin{prop}\label{p:HKE-generalu}
 For any  $T>0$, 
 there exists a constant $C>0$  depending on $T$ such that, for all 
$t\in (0, T]$   and $x,y \in D$,
	\begin{align*}
		&	p^\kappa(t,x,y)\\
		&\le C  \left(1 \wedge \frac{\delta_D(x)}{t^{1/\alpha}}
		\right)^p
		\left(1 \wedge \frac{\delta_D(y)}{t^{1/\alpha}}
		\right)^p
		\left[t^{-d/\alpha} \wedge  \left( \frac{ t }{|x-y|^{d+\alpha}}\Phi\left( \frac{(|x-y| \wedge A_0)^2}{
			(\delta_D(x,t) \wedge A_0)(\delta_D(y,t) \wedge A_0) }\right) \right)\right].
	\end{align*}	
\end{prop}
\begin{proof}
We assume, without loss of generality that $T\ge \wh R^\alpha$. 
 Let $t\in (0,T]$ and $x,y \in D$. 
 Let $n\in \N$ be such that 
	  $(n-1) \wh R^\alpha \le t<n \wh R^\alpha$ and set $h:=t/(n+1)$. Clearly,   $n\le 1+  (T/\wh R^\alpha) $. Using the semigroup properties  of $p^\kappa$ and $\overline p$, 
	 Proposition \ref{p:blow-up-step2} and Theorem \ref{t:CKSV25} (twice), we obtain that
	 	\begin{align*}
	 & p^\kappa(t,x,y)=
	 \int_D \int_D
	  p^\kappa(h/4,x,z) p^\kappa(t- h/2, z,w)  p^\kappa(h/4,w,y) \,dwdz\\
	&  \le c_1 \left(1 \wedge \frac{\delta_D(x)}{(h/4)^{1/\alpha}}
	\right)^p 	\left(1 \wedge \frac{\delta_D(y)}{(h/4)^{1/\alpha}}
		\right)^p
  \int_D \int_D
 \wt p(h/4,x,z) \overline p(t- h/2, z,w)  \wt p(h/4,w,y) \,dwdz\\
& \le 	c_2 n^{2p/\alpha} \left(1 \wedge \frac{\delta_D(x)}{t^{1/\alpha}}
	\right)^p 	\left(1 \wedge \frac{\delta_D(y)}{t^{1/\alpha}}
		\right)^p
 \int_D \int_D
\overline p(h/4,x,z) \overline p(t- h/2, z,w)  \overline p(h/4,w,y) \,dwdz\\
	& = 	c_2n^{2p/\alpha} \left(1 \wedge \frac{\delta_D(x)}{t^{1/\alpha}}
	\right)^p 	\left(1 \wedge \frac{\delta_D(y)}{t^{1/\alpha}}
		\right)^p	\overline p(t,x,y) \le c_3 \left(1 \wedge \frac{\delta_D(x)}{t^{1/\alpha}}
		\right)^p 	\left(1 \wedge \frac{\delta_D(y)}{t^{1/\alpha}}
		\right)^p \wt p (t, x,y),
	\end{align*}
	for some $c_3>0$ depending on $T$ but independent of $t,x,y$. The proof is complete.
\end{proof}


\section{Small time heat kernel lower bound}\label{s:lHKE}

The following proposition is the main result of this section.

\begin{prop}\label{p:lHKE-small}
 	Let $x_0 \in \overline D$ and $R>0$.   For any  $T>0$, 
there exists a constant $C>0$  depending on $T$,  but independent of $x_0$ and $R$, such that for all
$t\in (0, T\wedge R^\alpha]$   and $x,y \in B_D(x_0,R/2)$, \bk  	
\begin{align*}
		p^\kappa_{B_D(x_0,R)}(t,x,y) &\ge C \left(1 \wedge \frac{\delta_D(x)}{t^{1/\alpha}}\right)^p		\left(1 \wedge \frac{\delta_D(y)}{t^{1/\alpha}} \right)^p\nn\\		&\quad \times \left[t^{-d/\alpha} \wedge  \left( \frac{ t }{|x-y|^{d+\alpha}}	\Phi\left( \frac{(|x-y| \wedge A_0)^2}	{(\delta_D(x,t) \wedge A_0)(\delta_D(y,t) \wedge A_0) }\right) \right)\right].	
\end{align*}	
\end{prop}

The proof of Proposition \ref{p:lHKE-small} will be given at the end of this section.

Throughout this section, we fix $x_0 \in \overline D$ and $R>0$. We write $B:=B_D(x_0,R)$ for simplicity. 
{\it All constants appearing in the statements below will be chosen independently  of $x_0$ and $R$.}

\smallskip

 To prove Proposition \ref{p:lHKE-small}, we first recall the notion of a lifted point introduced in  \cite[Section 5.1]{CKSV25b}.
 For $x\in D$, we fix a point $Q_x\in \partial D$  satisfying $|x-Q_x|=\delta_D(x)$. Define $\mathbf{n}_x:=(x-Q_x)/|x-Q_x|$. Since $D$ is $C^{1,\eta}$ (hence Lipschitz), there exist  constants $\wh r \in (0, \wh R]$  and $ \eta_0 \in (0,1)$  such that\begin{equation}\label{e:lifting-property}	 \eta_0 (\delta_D(x) + r) \le 	\delta_D(x + r \bn_x)  \le \delta_D(x) + r \quad 	 \text{for all $x\in D$ and $r\in (0, \wh r ]$}.\end{equation} In particular, we have $x+r\bn_x\in D$ for all $x\in D$ and $r\in (0, \wh r]$. 

For $x\in D$ and $t>0$, the  \textit{lifted point} is defined as \begin{equation}\label{e:def-x(t)}	x(t):= x+ (t^{1/\alpha}/2)\bn_x. \end{equation} Note that, by \eqref{e:lifting-property}, we have  $x(t)\in D$ for all $x\in D$ and $t\in (0, (2 \wh r )^\alpha]$. Moreover, for all  $x\in D$ and $t\in (0, \wh r^\alpha]$, it holds that\begin{align}\label{e:lifting-standard}	\delta_D(y)  >2^{-2}\eta_0 (t/4)^{1/\alpha} \quad \text{for all  $y \in B(x(t/4), 2^{-2}  \eta_0 (t/4)^{1/\alpha})$.} \end{align}

We begin the proof of Proposition \ref{p:lHKE-small} with an interior near-diagonal lower bound for $p^\kappa$, which  follows from Proposition \ref{p:ndl2} and the semigroup property.

\begin{lemma}\label{l:ndl-refined} For any $\eps\in (0,1]$, there exists $C=C(\eps)>0$ such that  
for all $t\in (0,  (\wh r\wedge (R/2))^\alpha]$ and  
	$x\in B_D(x_0,3R/4)$ with $\delta_D(x) \ge \eps t^{1/\alpha}$,
	\begin{equation}\label{e:ndl-refined-claim}	
\inf_{y \in B(x, (\delta_D(x) \wedge t^{1/\alpha})/8) } 		p^\kappa_B(t,x,y)  \ge  C t^{-d/\alpha}.
	\end{equation}
\end{lemma}
\begin{proof} 
Let  $t\in (0, (\wh r\wedge (R/2))^\alpha]$ and  $x \in B_D(x_0,3R/4)$ with $\delta_D(x) \ge \eps t^{1/\alpha}$. 
Set $t_0:=2^{-2\alpha} (\delta_D(x) \wedge t^{1/\alpha})^\alpha$.  
Applying Proposition \ref{p:ndl2} 
with $\wh r\wedge (R/2)$ instead of $R$, $t=s$, $r=2s^{1/\alpha}$ and $b=1/2$, 
 we get that for all $s\in (0, t_0]$,  
	\begin{align}\label{e:ndl-refined-1}
	\inf_{y,z\in B(x,s^{1/\alpha}/2)}	p^\kappa_B(s,z,y) \ge 	\inf_{y,z\in B(x,s^{1/\alpha}/2)}	p^{\kappa}_{ B(x, 
		2	s^{1/\alpha}
			)}(s,z,y) \ge  c_1s^{-d/\alpha}.
	\end{align}
Here we use 
$R\ge 2 (\wh r\wedge (R/2))$ so that $B(x,2s^{1/\alpha}) \subset B$.

 Let $n\ge 1$ be such that  $nt_0 < t \le (n+1)t_0$.  	Note that  $n<t/t_0 \le  2^{2\alpha}\eps^{-\alpha}=:c_2(\eps)$.  Set $h:=(t-t_0)/n \in (0,t_0]$. 
	 Applying \eqref{e:ndl-refined-1} with $s= h$ and  $s= t_0$, we obtain
\begin{align*}
& 	\inf_{z,z' \in B(x, h^{1/\alpha}/2)}  p^\kappa_B(h, z, z' ) \ge  c_1 h^{-d/\alpha},\\
&	
	\inf_{z \in B(x, h^{1/\alpha}/2),\, y \in B(x, t_0^{1/\alpha}/2)}  p^\kappa_B(t_0, z, y ) \ge  c_1 t_0^{-d/\alpha} \ge  2^{-2d}c_1 t^{-d/\alpha}.
\end{align*}
Hence, using the semigroup property, we arrive at
	\begin{align*}
		&	p^\kappa_B(t,x,y) \ge \int_{B(x,h^{1/\alpha}/2)\times \cdots B(x,h^{1/\alpha}/2)} p^\kappa_B(h,x,z_1) \cdots p^\kappa_B(h,z_{n-1},z_{n}) \, p^\kappa_B(t_0, z_{n}, y )\,  dz_1\cdots dz_{n}\\
		&\ge (c_1h^{-d/\alpha})^{n} (2^{-2d}c_1t^{-d/\alpha})  \, m_d\big(B(x,h^{1/\alpha}/2)\big)^{n} =c_3t^{-d/\alpha} e^{-c_4(n-1)} \ge c_3t^{-d/\alpha} e^{-c_4 (c_2(\eps)-1)}. 
	\end{align*}
	The proof is complete.
\end{proof}

Combining Lemma \ref{l:ndl-refined} with the strong Markov property and the L\'evy system formula \eqref{e:Levysystem-X-kappa}, we obtain the following full lower heat kernel bound between neighborhoods of lifted points.

\begin{lemma}\label{l:lift-interior-estimate}
	Let $t\in (0,(\wh r\wedge (R/2))^\alpha]$ and   $x,y \in B_D(x_0,R/2)$. 
	There exists $C>0$ independent of $t,x,y$ such that
	\begin{align*}
		p^\kappa_B(t/2, z,w) \ge  C\left[t^{-d/\alpha} \wedge  \left( \frac{ t }{|x-y|^{d+\alpha}}
		\Phi\left( \frac{(|x-y| \wedge A_0)^2}{(\delta_D(x,t) \wedge A_0)(\delta_D(y,t) \wedge A_0) }\right) \right)\right]
	\end{align*}
	for all	$z\in B(x(t/4), 2^{-6} \eta_0 (t/8)^{1/\alpha})$ and 	$w\in B(y(t/4), 2^{-6}  \eta_0 (t/8)^{1/\alpha})$.
\end{lemma}

\begin{proof} 
Without loss of generality, we  assume $\delta_D(x)\le \delta_D(y)$. Write $r:=2^{-6} \eta_0(t/8)^{1/\alpha}$. Let  $z\in B(x(t/4),r)$ and 	$w\in B(y(t/4), r)$. 
Note that $z,w \in B_D(x_0,3R/4)$.
Recall from \eqref{e:lifting-standard} that  $\delta_D(z) \wedge \delta_D(w) > 2^{-2}\eta_0 (t/4)^{1/\alpha}$, which implies that, for
$s\in [t/6,t/2]$, we have 
$\delta_D(z) \wedge \delta_D(w) \wedge s^{1/\alpha} > 2^{-2-1/\alpha}\eta_0 s^{1/\alpha}$ and 
$\delta_D(z) \wedge \delta_D(w) \wedge s^{1/\alpha}>2^4r$.
Thus,  by Lemma \ref{l:ndl-refined}, it holds that 
\begin{align}\label{e:lift-interior-estimate-1}
	\bigg(\inf_{s\in [t/6,t/2],\, |z-z'|<2r}	p^\kappa_B(s, z, z') \bigg) \wedge  \bigg( \inf_{s\in [t/6,t/2],\, |w-w'|<2r}	p^\kappa_B(s,  w,w') \bigg) \ge c_1 t^{-d/\alpha}.
\end{align}
If $|z-w|<2r$, then the result follows from \eqref{e:lift-interior-estimate-1}.

Suppose $|z-w|\ge 2r$. For all $z'\in B(z,r/2)$, by \eqref{e:lifting-property}, 
$$
\delta_D(z') \ge \delta_D(x(t/4)) -3r/2 > \eta_0 
\delta_D(x) + 2^{-2}\eta_0 
(t/4)^{1/\alpha} \ge 2^{-2-2/\alpha} \eta_0 \delta_D(x,t)
$$
and  
$\delta_D(z') \le \delta_D(x(t/4)) + 3r/2  \le 2\delta_D(x,t)$. Similarly,  we have  $2^{-2-2/\alpha}\eta_0 \delta_D(y,t)\le \delta_D(w') \le  2\delta_D(y,t)$ for all 
$w'\in B(w,r/2)$. 
Further, for all $z'\in B(z,r/2)$ and $w'\in B(w,r/2)$, since
\begin{align*}
	|z'-w'| < |x(t/4)-y(t/4)|+3r\le  |x-y| + |x-x(t/4)|+ |y-y(t/4)| +3r \le |x-y|+c_2t^{1/\alpha}
\end{align*}
and 
$$|z'-w'| >(|z-w|-r) \vee (|x-y| - |x-x(t/4)| -  |y-y(t/4)| -3r)  \ge  r \vee      (|x-y|- c_3t^{1/\alpha}),$$ 
we get that
\begin{align*}
	|z'-w'|  \asymp \begin{cases}
		t^{1/\alpha} &\mbox{ if $|x-y|\le 2c_3t^{1/\alpha}$},\\
		|x-y| &\mbox{ if $|x-y|> 2c_3t^{1/\alpha}$}
	\end{cases}.
\end{align*}
Thus, for all  $z'\in B(z,r/2)$ and $w'\in B(w,r/2)$,  if $|x-y| \le 2c_3t^{1/\alpha}$, then
\begin{align}\label{e:lift-interior-estimate-2}
	\frac{\sB(z',w')}{|z'-w'|^{d+\alpha}} &\ge \frac{c_4}{t^{(d+\alpha)/\alpha}} 
	\Phi\,\bigg( \frac{(t^{1/\alpha} \wedge A_0)^2}{(\delta_D(x,t) \wedge A_0)(\delta_D(y,t) \wedge A_0) }\bigg)
	\nn \\ 
	&\ge  \frac{c_5}{t^{(d+\alpha)/\alpha}}\Phi\left( \frac{(|x-y| \wedge A_0)^2}{(\delta_D(x,t) \wedge A_0)(\delta_D(y,t) \wedge A_0) }\right).
\end{align}
If $|x-y| > 2c_3t^{1/\alpha}$, then
\begin{align}\label{e:lift-interior-estimate-3}
	\frac{\sB(z',w')}{|z'-w'|^{d+\alpha}} 
	\ge  \frac{ c_6}{|x-y|^{d+\alpha}}\Phi\left( \frac{(|x-y| \wedge A_0)^2}{(\delta_D(x,t) \wedge A_0)(\delta_D(y,t) \wedge A_0) }\right).
\end{align}
Using the strong Markov property in the first line below, \eqref{e:lift-interior-estimate-1} in the second,  the L\'evy system formula \eqref{e:Levysystem-X-kappa} in the third,  and  \eqref{e:lift-interior-estimate-2} and \eqref{e:lift-interior-estimate-3} in the fourth, we deduce that 
\begin{align}\label{e:lift-interior-estimate-4}
	& p^\kappa_B(t/2,z,w) \ge \E_z \bigg[ p^\kappa_B(t/2- 
	\tau^\kappa_{B(z,r/2)}, X^{\kappa,B}_{\tau^\kappa_{B(z,r/2)}}, w);\,  \tau^\kappa_{B(z,r/2)} \le t/3, \, 
		X^{\kappa,B}_{\tau^\kappa_{B(z,r/2)}} \in B(w,r/2)  \bigg]\nn \\
	&\ge c_1t^{-d/\alpha} \P_z \Big( \tau^\kappa_{B(z,r/2)} \le t/3, \, X^{\kappa}_{\tau^\kappa_{B(z,r/2)} \wedge \tau^\kappa_B} \in B(w,r/2)  \Big)\nn\\
	&\ge c_1t^{-d/\alpha} \E_z \bigg[ \int_{0}^{	\tau^\kappa_{B(z,r/2)} 
	\wedge (t/3)} \int_{B(w,r/2)}\frac{\sB(X^{\kappa}_s, w')}{|X^{\kappa}_s-w'|^{d+\alpha}} dw\,ds \bigg]\nn\\
	&\ge \frac{c_7 }{(|x-y| \vee t^{1/\alpha})^{d+\alpha}} \E_z [	\tau^\kappa_{B(z,r/2)} 
	\wedge (t/3)] \Phi\left( \frac{(|x-y| \wedge A_0)^2}{(\delta_D(x,t) \wedge A_0)(\delta_D(y,t) \wedge A_0) }\right).
\end{align}
By applying Proposition \ref{p:ndl2} with $b=1/2$, we get that 
\begin{align}\label{e:lift-interior-estimate-5}
	&	\E_z [ \tau^\kappa_{B(z,r/2)} \wedge (t/3)] \ge (r/4)^\alpha \P_z(\tau^\kappa_{B(z,r/2)}
        \ge (r/4)^\alpha ) \ge (r/4)^\alpha \P_z\big(X^{\kappa, B(z,r/2)}_{(r/4)^\alpha} \in B(z, r/8) \big) \nn\\
	&= (r/4)^\alpha \int_{B(z,r/8)} p^{\kappa}_{ B(z,r/2)}((r/4)^\alpha, z, z')dz' 
	\ge c_8 (r/4)^{\alpha-d}\int_{B(z,r/8)}dz'  = c_9r^\alpha = c_{10}t.
\end{align}
Combining \eqref{e:lift-interior-estimate-4} and  \eqref{e:lift-interior-estimate-5} with the fact that 
$t/(|x-y|\vee t^{1/\alpha})^{d+\alpha} = t^{-d/\alpha} \wedge (t|x-y|^{-d-\alpha})$, we obtain the desired result. 
\end{proof}

The following result is the counterpart of \cite[Lemma 5.4]{CKSV25b} in the current setting.

\begin{lemma}\label{l:lift-probability}
There exist constants 
$\lambda \in (0,2^{-2-3/\alpha}]$
 and  $C>0$ such that for all $x\in B_D(x_0,R/2)$ and  
$t\in (0, (\wh r\wedge (R/2))^\alpha  ]$,
\begin{align*}
	\P_x \Big( X^{\kappa,B}_{\tau^\kappa_{B_{ D}(x, \lambda t^{1/\alpha})} \wedge (t/8)} 
	\in B(x(t/4), 2^{-6} \eta_0 (t/8)^{1/\alpha}) \Big)  
	\ge C\left( 1 \wedge \frac{\delta_D(x)}{t^{1/\alpha}}\right)^p.
\end{align*}	
\end{lemma}
\begin{proof}
Let  $x\in B_D(x_0,R/2)$, $t\in (0, (\wh r\wedge (R/2))^\alpha]$, and let 
$\lambda\in (0,2^{-2-3/\alpha} ]$
 be a constant to be chosen later.   
Note that  $2^{-2-3/\alpha} t^{1/\alpha}  + 2^{-6} \eta_0 (t/8)^{1/\alpha} <(17/64)(t/4)^{1/\alpha} =(17/32) |x-x(t/4)| $. 
Since $\lambda \le 2^{-2-3/\alpha}$, it follows that
\begin{align*}
		|z-y| \asymp t^{1/\alpha}  \quad \text{for  $z\in B_{ D}(x, \lambda t^{1/\alpha})$ and $y \in B(x(t/4), 2^{-6}
		  \eta_0 
		(t/8)^{1/\alpha})$.} 
\end{align*}
Hence,   using the fact that $\Phi \ge 1$, we get from \eqref{e:J-Phi} that for all  $z\in B_{ D}(x, \lambda t^{1/\alpha})$ and $y \in B(x(t/4), 2^{-6}
\eta_0 
(t/8)^{1/\alpha})$,
\begin{align}\label{e:lift-probability-1}
	\frac{\sB(z,y)}{|z-y|^{d+\alpha}} \ge \frac{C_1^{-1}}{|z-y|^{d+\alpha}} \ge c_2 t^{-d/\alpha-1}.
\end{align}
Combining the L\'evy system formula \eqref{e:Levysystem-X-kappa} with \eqref{e:lift-probability-1}, since $\tau^\kappa_{B_{ D}(x, \lambda t^{1/\alpha})}\le \tau^\kappa_B$,  we obtain
\begin{align}\label{e:lift-probability-2}
	&\P_x \Big( X^{\kappa,B}_{ \tau^\kappa_{B_{ D}(x, \lambda t^{1/\alpha})} \wedge (t/8)} \in B(x(t/4), 2^{-6}\eta_0 (t/8)^{1/\alpha})  \Big)\nn\\
	&= \E_x \bigg[\int_0^{ \tau^\kappa_B \wedge \tau^\kappa_{B_{ D}(x, \lambda t^{1/\alpha})} \wedge (t/8) }\int_{B(x(t/4), 2^{-6}\eta_0(t/8)^{1/\alpha}) } 
		 \frac{\sB(X^\kappa_s, y)}{|X^\kappa_s-y|^{d+\alpha}}\,dy ds\bigg] \nn\\
    &\ge  c_2 t^{-d/\alpha-1}  m_d\big(B(x(t/4), 2^{-6}	    \eta_0 (t/8)^{1/\alpha}) \big) \E_x \big[ \tau^\kappa_{B_{ D}(x, \lambda t^{1/\alpha})} \wedge (t/8) \big]\nn\\
	& \ge c_3 t^{-1}\E_x \big[ \tau^\kappa_{B_{ D}(x, \lambda t^{1/\alpha})};\tau^\kappa_{B_{ D}(x, \lambda t^{1/\alpha})}\le t/8 \big] \nn\\
	&= c_3 t^{-1}\E_x \big[\tau^\kappa_{B_{ D}(x, \lambda t^{1/\alpha})} \big] - c_3t^{-1} 
		\E_x \big[ \tau^\kappa_{B_{ D}(x, \lambda t^{1/\alpha})} ; \tau^\kappa_{B_{ D}(x, \lambda t^{1/\alpha})} > t/8 \big].
\end{align}  
  Using Lemma \ref{l:survival-exponential}  and the fact that $ \sup_{a\ge 0}a^ke^{-a} = (k/e)^k$ for all $k>0$, we obtain
\begin{align}\label{e:lift-probability-3}
	&t^{-1}\E_x \big[ \tau^\kappa_{B_{ D}(x, \lambda t^{1/\alpha})}  ;\, \tau^\kappa_{B_{ D}(x, \lambda t^{1/\alpha})} > t/8 \big]\nn\\
	& \le \sum_{n=1}^\infty \frac{n}{4}\P_x \big( \tau^\kappa_{B_{ D}(x, \lambda t^{1/\alpha})}
		\in (nt/8, (n+1)t/8]\big)\le c_4  \left( 1 \wedge \frac{\delta_D(x)}{\lambda t^{1/\alpha}}\right)^p
		\sum_{n=1}^\infty n e^{-c_5n/\lambda^\alpha} \nn\\
	&\le c_6 \lambda^{3\alpha} \left( 1 \wedge \frac{\delta_D(x)}{\lambda t^{1/\alpha}} \right)^p
		\sum_{n=1}^\infty n^{-2} = c_7 \lambda^{3\alpha} \left( 1 \wedge \frac{\delta_D(x)}{\lambda t^{1/\alpha}}\right)^p, 
\end{align} 
where $c_7$ is independent of $\lambda$. Combining \eqref{e:lift-probability-2} with Proposition \ref{p:met-bound-new}  and  \eqref{e:lift-probability-3}, we arrive at
\begin{align*}
	\P_x \Big( X^{\kappa,B}_{ \tau^\kappa_{B_{ D}(x, \lambda t^{1/\alpha})} \wedge (t/8)} \in B(x(t/4), 2^{-6}\eta_0 (t/8)^{1/\alpha})  \Big) \ge  (c_8\lambda^{\alpha}-c_7 \lambda^{3\alpha})\left( 1 \wedge \frac{\delta_D(x)}{\lambda t^{1/\alpha}}
	\right)^p.
\end{align*}
By choosing $\lambda$ sufficiently small, this yields the desired result.
\end{proof}

The following result is the counterpart of \cite[Lemma 5.5]{CKSV25b} in the current setting.

\begin{lemma}\label{l:lift}
There exists $C>0$ such that  for all  $x\in B_D(x_0,R/2)$ and  $t\in (0,  (\wh r\wedge (R/2))^\alpha ]$, 
$$ 	
\inf_{y \in B(x(t/4), \,2^{-6}  \eta_0 (t/8)^{1/\alpha})} 		p^\kappa_B(t/4,x,y) \ge C  \left(1 \wedge \frac{\delta_D(x)}{t^{1/\alpha}}\right)^p t^{-d/\alpha}.
$$
\end{lemma}
\begin{proof}
Let $x\in B_D(x_0,R/2)$,  $t\in (0,  (\wh r\wedge (R/2))^\alpha ]$ and $y \in B(x(t/4), \,2^{-6}  \eta_0 (t/8)^{1/\alpha})$.
By \eqref{e:lifting-standard}, 
for all $s\in [t/8,t/4]$ and $z\in B(x(t/4),2^{-6}\eta_0 (t/8)^{1/\alpha})$,  we have
$\delta_D(y) > 2^{-2}\eta_0 (t/4)^{1/\alpha} \ge 2^{-2} \eta_0 s^{1/\alpha}$ and 
$$
|y-z| \le |y- x(t/4)|+|x(t/4)-z| <2^{-5}\eta_0 
(t/8)^{1/\alpha}  <2^{-3}(\delta_D(y) \wedge s^{1/\alpha}).
$$
Thus,  we can apply Lemma \ref{l:ndl-refined} with $\eps=2^{-2}\eta_0 $ to obtain 
\begin{align}\label{e:lift-on-diagonal-1}
	p^\kappa_B(s,z,y) \ge  c_1 t^{-d/\alpha} \quad \text{
		for all $s\in [t/8,t/4]$ and $z\in B(x(t/4),2^{-6} \eta_0 (t/8)^{1/\alpha})$.  } 
\end{align}
Let $\lambda$ be the constant in Lemma \ref{l:lift-probability} and write $V:=B_D(x,\lambda t^{1/\alpha})$.
Using the strong Markov property, \eqref{e:lift-on-diagonal-1} and Lemma \ref{l:lift-probability}, we obtain
\begin{align*}
&	p^\kappa_B(t/4,x,y) \ge 	\E_x\left[p^\kappa_B(t/4-\tau^\kappa_{V}, X^{\kappa,B}_{\tau^\kappa_V}, y);\,\tau^\kappa_{V} \le t/8,\, X^{\kappa,B}_{\tau^\kappa_V} 		\in B(x(t/4), 2^{-6} \eta_0 (t/8)^{1/\alpha})\right]\nn\\	&\ge \bigg(\inf_{  s\in [t/8,t/4], \, z \in B(x(t/4), 2^{-6} \eta_0 (t/8)^{1/\alpha})} p^\kappa_B (s, z,y) 	\bigg)		\P_x \Big( X^{\kappa,B}_{\tau^\kappa_{V} \wedge (t/8)} \in B(x(t/4), 2^{-6} \eta_0 (t/8)^{1/\alpha}) \Big)  \nn\\	&\ge c_2 \left(1 \wedge \frac{\delta_D(x)}{t^{1/\alpha}}\right)^pt^{-d/\alpha}.
\end{align*}
\end{proof}

We will use the following elementary inequality:
\begin{align}\label{e:simple}
	\frac{T \wedge R^\alpha}{\wh r^\alpha \wedge (R/2)^\alpha} \le \frac{T}{\wh r^\alpha} \vee 2^\alpha \quad \text{for  
	all $R, T>0$}.
\end{align}

Recall that $\wt p(t,x,y)$ is defined in \eqref{e:def-wtp}.  
\begin{proof}[Proof of Proposition \ref{p:lHKE-small}]
Let $t\in (0, T \wedge R^\alpha]$ and  $x,y \in B_D(x_0,R/2)$.  
We assume, without loss of generality that $T\ge \wh r^\alpha$, and write  $r_0:=\wh r \wedge (R/2)$  for simplicity.  Recall that $B=B_D(x_0,R)$.  

Define $V(z,s):=B(z(s/4), 2^{-6}\eta_0 (s/8)^{1/\alpha})$ for $z\in D$ and  $s\in (0, \wh r^\alpha]$. 
Let $n\in \N$ be such that $(n-1) r_0^\alpha \le t<n r_0^\alpha$ and set $h:=t/(n+1)$. Note that 
$n-1 \le t/r_0^\alpha \le   (T/\wh r^\alpha) \vee 2^\alpha$ by \eqref{e:simple}.  Thus,  
using \eqref{e:Phi-scaling} and the fact that $n\le 1+(T/\wh r^\alpha) \vee 2^\alpha$, we get that 
there are comparison constants depending on $T$, independent of $t$, such that 
\begin{align}\label{e:lHKE-1}
	 \wt  p(h, z,w)\asymp \wt  p(t, z,w) \quad \text{for all $z,w \in D$}.
\end{align}
By Lemmas \ref{l:lift} and \ref{l:lift-interior-estimate} and the symmetry of $p^\kappa_B$, we have
\begin{align*}
\inf_{z\in V(x,h)}	p^\kappa_B(h/4, x,z) &\ge c_1 \left(1 \wedge \frac{\delta_D(x)}{h^{1/\alpha}}\right)^p h^{-d/\alpha},\\
 \inf_{w\in V(y,h)}	p^\kappa_B(h/4, w,y) &\ge c_1 \left(1 \wedge \frac{\delta_D(y)}{h^{1/\alpha}}\right)^p h^{-d/\alpha},\\
	 \inf_{z,z'\in V(x,h)}	p^\kappa_B(h/2, z,z')&\ge c_2 h^{-d/\alpha},\quad \;\;
	  \inf_{z\in V(x,h), \, w \in V(y,h)}	p^\kappa_B(h/2, z,z') \ge c_3 \wt p(h,x,y).
\end{align*}
Hence, using the semigroup property, we obtain that  
	\begin{align}\label{e:hkeln1}
&		p^\kappa_B(t,x,y)\ge \int_{V(x,h)} \cdots \int_{V(x,h)} \int_{V(y,h)} p^\kappa_B(h/4,x,z_0) \bigg(\prod_{j=1}^{2n}p^\kappa_B(h/2,z_{j-1},z_j)  \bigg) \nn \\
		&\qquad\qquad\qquad\qquad\qquad\qquad \times p^\kappa_B(h/2,z_{2n},w) p^\kappa_B(h/4,w,y) \,dw\,
		dz_{2n} \cdots dz_{0}
		\nn\\	
		&\ge c_1^2 c_2^{2n} c_3 m_d(V(x,h))^{2n+1} m_d(V(y,h)) h^{-(2n+2)d/\alpha} \left(1 \wedge \frac{\delta_D(x)}{h^{1/\alpha}}\right)^p\left(1 \wedge \frac{\delta_D(y)}{h^{1/\alpha}}\right)^p \wt p(h,x,y)\nn \\
		&= c_1^2 c_2^{2n} c_5^{2n+2} \left(1 \wedge \frac{\delta_D(x)}{h^{1/\alpha}}\right)^p\left(1 \wedge \frac{\delta_D(y)}{h^{1/\alpha}}\right)^p \wt p(h,x,y).
	\end{align} 
Using  $n\le 1+(T/\wh r^\alpha) \vee 2^\alpha$, we get
$$c_1^2 c_2^{2n} c_5^{2n+2} \left(1 \wedge \frac{\delta_D(x)}{h^{1/\alpha}}\right)^p\left(1 \wedge \frac{\delta_D(y)}{h^{1/\alpha}}\right)^p \le c_6 \left(1 \wedge \frac{\delta_D(x)}{t^{1/\alpha}}\right)^p\left(1 \wedge \frac{\delta_D(y)}{t^{1/\alpha}}\right)^p
$$
for some $c_6>0$ depending on $T$ but independent of $t,x,y, R$. Combining this with \eqref{e:hkeln1} and \eqref{e:lHKE-1} yields the desired result.
\end{proof}


\section{Proofs of main results and Green function estimates }\label{s:main-proofs}

Combining the lower heat kernel bounds established in Section  \ref{s:lHKE} with Lemma  \ref{l:life}, we obtain the following sharp two-sided estimates for the lifetime.
\begin{prop}\label{p:life}	For any $T>0$,  there exists  $C=C(T) >1$ such that$$ C^{-1} \left(1 \wedge \frac{\delta_D(x)}{t^{1/\alpha}}\right)^p\le \P_x (\zeta^\kappa>t) \le  C \left(1 \wedge \frac{\delta_D(x)}{t^{1/\alpha}}\right)^p \quad \text{for all $t\in (0,T]$ and $x\in D$}.$$\end{prop}
\begin{proof}
	The upper bound was established in Lemma \ref{l:life}.  We prove the lower bound. 	Let $t\in (0,T]$, $x\in D$ and let $c_1$ denote the constant $C$ in Proposition \ref{p:met-bound-new}.
	Set $r_1:=\wh r \wedge (c_1^{1/\alpha}\wh R)$. 			If $t< r_1^\alpha$, then the result follows from Proposition \ref{p:met-bound-new} since  
	$\zeta^\kappa \ge \tau^\kappa_{B_D (x,r)}$. 
	Suppose $r_1^\alpha \le t\le T$. Set 
	$V:=B(x(r_1^\alpha/8), \,2^{-6}  \eta_0 (r_1^\alpha/16)^{1/\alpha})$ where $ \eta_0  \in (0, 1)$ is the constant in \eqref{e:lifting-property}.
	Using  Lemma \ref{l:lift},  
	Proposition \ref{p:lHKE-small} with the fact that $\Phi \ge 1$ and  \eqref{e:lifting-standard}, we obtain that for all $z,y \in V$,
	\begin{align}
		p^\kappa_{B_D(x,T^{1/\alpha})}( r_1^\alpha/2,x,z)  &\ge c_2 \left(1 \wedge \frac{\delta_D(x)}{r_1}\right)^p r_1^{-d} \ge  c_2 \left(1 \wedge \frac{\delta_D(x)}{t^{1/\alpha}}\right)^p r_1^{-d} ,\\
		p^\kappa_{B_D(x,T^{1/\alpha})}( t- r_1^\alpha/2,z,y) &\ge c_3 \left(1 \wedge \frac{\delta_D(z)}{T^{1/\alpha}}\right)^p\left(1 \wedge \frac{\delta_D(y)}{T^{1/\alpha}}\right)^p \left( T^{-d/\alpha} \wedge \frac{r_1^\alpha/2}{|z-y|^{d+\alpha}} \right) \ge c_4.
	\end{align}
	Thus, using the semigroup property, we get that
	\begin{align*}	
		&\P_x (\zeta^\kappa>t) \ge \P_x (\tau^\kappa_{B_D(x,T^{1/\alpha})}>t)  \ge \int_{V} p^\kappa_{B_D(x,T^{1/\alpha})}( t,x,y)  dy  \\
		&\ge  \int_{V} \int_{V}  p^\kappa_{B_D(x,T^{1/\alpha})}( r_1^\alpha/2,x,z)   p^\kappa_{B_D(x,T^{1/\alpha})}( t-r_1^\alpha/2,z,y)     dy dz \ge c_5 m_d(V)^2 \left(1 \wedge \frac{\delta_D(x)}{t^{1/\alpha}}\right)^p.	\end{align*}
	The proof is complete.
\end{proof}

Combining Proposition \ref{p:HKE-generalu}, Proposition \ref{p:lHKE-small} and the monotonicity of heat kernels
we  obtain the following slight generalization  of Theorem \ref{t:HKE}. 

\begin{thm}\label{t:HKE-general}
 Let $x_0 \in \overline D$ and $R>0$.   For any  $T>0$, there exist comparison constants depending on $T$,  but independent of $x_0$ and $R$ such that, for all 
$t\in (0, T\wedge R^\alpha]$   and $x,y \in B_D(x_0,R/2)$,
	\begin{align*}
		&	p^\kappa_{B_D(x_0,R)}(t,x,y)\\
		&\asymp \left(1 \wedge \frac{\delta_D(x)}{t^{1/\alpha}}
		\right)^p
		\left(1 \wedge \frac{\delta_D(y)}{t^{1/\alpha}}
		\right)^p
		\left[t^{-d/\alpha} \wedge  \left( \frac{ t }{|x-y|^{d+\alpha}}\Phi\left( \frac{(|x-y| \wedge A_0)^2}{
			(\delta_D(x,t) \wedge A_0)(\delta_D(y,t) \wedge A_0) }\right) \right)\right].
	\end{align*}	
\end{thm}

\medskip
	
	The proof of Theorem \ref{t:HKE} is now immediate.
	
	\begin{proof}[Proof of Theorem \ref{t:HKE}.]
		Let $x_0 \in \overline D$.	By the monotone convergence theorem, $p^\kappa(t,x,y) = \lim_{R\to \infty} $ $p^\kappa_{B_D(x_0,R)}(t,x,y)$ for all $t>0$ and $x,y \in D$. Thus, the result follows from Theorem \ref{t:HKE-general} together with Proposition  \ref{p:life}.
	\end{proof}

The proof of Theorem \ref{t:main2} is standard and we only give a brief outline.

\begin{proof}[Proof of Theorem \ref{t:main2}.]
When $D$ is bounded, the proof of the large time estimates in \eqref{t:e:Lhke} using the estimates in Theorem \ref{t:HKE} is well known, see the proof of \cite[Theorem 1.7(ii)]{CKSV25b}. Since the proof is the  same as that of \cite[Theorem 1.7(ii)]{CKSV25b}, we omit it.
\end{proof}

For an open set $U\subset D$, 
 the Green function $G^\kappa_U$ of $X^{\kappa,U}$ is defined by
$$
G^\kappa_U(x, y):=\int^\infty_0p^\kappa_U(t, x, y)dt, \quad x,y \in U.
$$
Having both small and large time estimates of the heat kernel suffices to obtain sharp two-sided estimates of the Green functions. 
\begin{thm}\label{t:gfe}
For any $x_0 \in \overline D$ and $R>0$, there exist comparison constants independent of $x_0$ such that for all 
$x,y \in B_D(x_0,R/2)$,
\begin{align*}
G^\kappa_{B_D(x_0,R)}(x,y)&\asymp 
\displaystyle 
\left(\frac{\delta_D(x)}{|x-y|}  \wedge 1 \right)^p\left(\frac{\delta_D(y)}{|x-y|}  \wedge 1 \right)^p\\
&\qquad \times 
\begin{cases}
 \displaystyle\frac{1}{|x-y|^{d-\alpha}},
&\alpha<d;\\[8pt]
\displaystyle\log\left(e+ \frac{ (\delta_D(x) \vee \delta_D(y)) \wedge R}{|x-y|}\right),
&\alpha=1=d;\\[8pt]
\left( (\delta_D(x) \vee \delta_D(y) \vee |x-y|) \wedge R\right)^{\alpha-1},
&\alpha>1=d.
\end{cases}
\end{align*} 
\end{thm}
\begin{proof}
Let $x_0 \in \overline D$, $R>0$ and $x,y \in B_D(x_0,R/2)$.	
For the upper bound, when $d>\alpha$ or $\delta_D(x) \vee \delta_D(y) < 2R$,  the result follows from \cite[Proposition 7.4]{CKSV26a}. 
 Suppose $\alpha \ge 1=d$ and $\delta_D(x) \vee \delta_D(y) \ge 2R$. Note that $\delta_D(x) \wedge \delta_D(y) \ge \delta_D(x) \vee \delta_D(y) - |x-y| \ge R$ in this case.   
 Let $\beta:=(\ub + \alpha)/2 \in (\ub,\alpha)$.  By Theorem \ref{t:HKE}, \eqref{e:without-A0}, and \eqref{e:Phi-scaling}, we get that for all $t\in (0,R^\alpha]$,
\begin{align*} 
	p^\kappa_{B_D(x_0,R)}(t,x,y) & \le  c_1
	\begin{cases}
		t^{-1/\alpha} &\mbox{if $|x-y|\le  t^{1/\alpha}$},\\
		t|x-y|^{-1-\alpha} \Phi(|x-y|^2/t^{2/\alpha})	&     \mbox{if $|x-y|> t^{1/\alpha}$},
	\end{cases}\\
	& \le  c_2
	\begin{cases}
		t^{-1/\alpha} &\mbox{if $|x-y|\le  t^{1/\alpha}$},\\
		t^{1-2\beta/\alpha}|x-y|^{-1-\alpha +2\beta}	&     \mbox{if $|x-y|> t^{1/\alpha}$}.
	\end{cases}
\end{align*}
Further, by  \cite[Lemma 4.15]{CKSV26a}, we have $	p^\kappa_{B_D(x_0,R)}(t,x,y)\le c_3e^{-\mu t}$ for all $t\ge R^\alpha$,  with constants $c_3>0$ and $\mu>0$ independent of $x_0$. Thus, since $\beta<\alpha$, we obtain
\begin{align*}
	G^\kappa_{B_D(x_0,R)} &\le c_2\int_0^{|x-y|^\alpha} \frac{t^{1-2\beta/\alpha}}{|x-y|^{1+\alpha-2\beta}} dt  + c_2 \int_{|x-y|^\alpha}^{R^\alpha} t^{-1/\alpha}dt + c_3 \int_{R^\alpha}^\infty e^{-\mu t}dt \\
	&\le c_4 |x-y|^{\alpha-1} + c_5\left( \ind_{\alpha=1} \log(R/|x-y|) + \ind_{\alpha>1} R^{\alpha-1} \right) + c_6 e^{-\mu R^\alpha}.
\end{align*}This yields the desired upper bound. 

 For the lower bound, using Theorem \ref{t:HKE-general}  (with $T=R^\alpha$)  and the fact that $\Phi \ge 1$, 
 we obtain for all $x,y \in B_D(x_0,R/2)$,
	\begin{align*}
	&	G^\kappa_{B_D(x_0,R)}(x,y) \ge \int_{|x-y|^\alpha/2}^{|x-y|^\alpha} p^\kappa _{B_D(x_0,R)}(t,x,y)dt\\
	&  \ge c_1\left(\frac{\delta_D(x)}{|x-y|}  \wedge 1 \right)^p\left(\frac{\delta_D(y)}{|x-y|}  \wedge 1 \right)^p \int_{|x-y|^\alpha/2}^{|x-y|^\alpha} \left(t^{-d/\alpha} \wedge \frac{t}{|x-y|^{d+\alpha}}\right) dt\\
	&= c_1\left(\frac{\delta_D(x)}{|x-y|}  \wedge 1 \right)^p\left(\frac{\delta_D(y)}{|x-y|}  \wedge 1 \right)^p \int_{|x-y|^\alpha/2}^{|x-y|^\alpha} \frac{t}{|x-y|^{d+\alpha}} dt\\
	& = c_2\left(\frac{\delta_D(x)}{|x-y|}  \wedge 1 \right)^p\left(\frac{\delta_D(y)}{|x-y|}  \wedge 1 \right)^p  |x-y|^{\alpha-d}.
	\end{align*}
	This proves the desired lower bound for the cases $\alpha<d$, or $\alpha\ge 1=d$ and $\delta_D(x) \vee \delta_D(y) \le 2|x-y|$.
	
	Suppose $\alpha\ge 1=d$ and let $x,y\in B_D(x_0,R/2)$ be such that  $\delta_D(x) \vee \delta_D(y) > 2|x-y|$. Note that $\delta_D(x) \wedge \delta_D(y) \ge \delta_D(x) \vee \delta_D(y)  - |x-y| > (\delta_D(x) \vee \delta_D(y))/2 >|x-y|$. 
	Using Theorem \ref{t:HKE-general} and the fact that $\Phi \ge 1$ again, we get that 
	\begin{align*}
		&	G^\kappa_{B_D(x_0,R)}(x,y) \ge 
	 \int_{|x-y|^\alpha/2}^{(\delta_D(x) \vee \delta_D(y))^\alpha \wedge R^\alpha} p^\kappa _{B_D(x_0,R)}(t,x,y)dt\\
		&  \ge c_3 \int_{|x-y|^\alpha/2}^{(\delta_D(x) \vee \delta_D(y))^\alpha \wedge R^\alpha} \left( 1 \wedge \frac{\delta_D(x) \wedge \delta_D(y)}{t^{1/\alpha}}\right)^p \left(t^{-d/\alpha} \wedge \frac{t}{|x-y|^{d+\alpha}}\right) dt\\
		&\ge c_3 \left( 1 \wedge \frac{\delta_D(x) \wedge \delta_D(y)}{\delta_D(x) \vee \delta_D(y)}\right)^p \int_{|x-y|^\alpha/2}^{(\delta_D(x) \vee \delta_D(y))^\alpha\wedge R^\alpha}  t^{-d/\alpha}  dt\\
		&\ge 2^{-p}c_4 \begin{cases}
		\displaystyle\log\left(\frac{ 
		2 ((\delta_D(x) \vee \delta_D(y)) \wedge R)	}{|x-y|}\right)	&\mbox{ if $\alpha=1=d$},\\[8pt]
			((\delta_D(x) \vee \delta_D(y))\wedge R)^{\alpha-1}	 - 
		2^{-(\alpha-1)/\alpha}  |x-y|^{\alpha-1}
				&\mbox{ if $\alpha>1=d$}.
		\end{cases}
	\end{align*}
	Since 
	$(\delta_D(x) \vee \delta_D(y)) \wedge R >|x-y|$, this yields the desired lower bound and completes the proof.
	
\end{proof}

\begin{corollary}\label{c:gfe}
If $D\subset \R^d$, $d\ge 1$, is a bounded $C^{1,\eta}$ open set, then the Green function $G^\kappa(x,y)$ for $X^{\kappa}$ has the following estimates:  for all $x,y \in D$,
\begin{equation}\label{e:FGE}
	G^\kappa (x,y)\asymp 
	\displaystyle 
	\left(\frac{\delta_D(x)}{|x-y|}  \wedge 1 \right)^p\left(\frac{\delta_D(y)}{|x-y|}  \wedge 1 \right)^p
	\begin{cases}
		\displaystyle\frac{1}{|x-y|^{d-\alpha}},
		&\alpha<d;\\[8pt]
		\displaystyle\log\left(e+ \frac{ \delta_D(x) \vee \delta_D(y)}{|x-y|}\right),
		&\alpha=1=d;\\[8pt]
		\left(\delta_D(x) \vee \delta_D(y) \vee |x-y|\right)^{\alpha-1},
		&\alpha>1=d.
	\end{cases}
\end{equation}  
\end{corollary}
\begin{proof}
The result follows directly from Theorem \ref{t:gfe} by taking $R=3\,\diam(D)$.
\end{proof}

\bigskip
\noindent
\textbf{Declarations:} On behalf of all authors, the corresponding author states that there is no conflict of
interest. Data sharing is not applicable to this article as no datasets were generated or analyzed during the
current study.

\vspace{.1in}
\begin{singlespace}

\end{singlespace}

\end{document}